\documentclass[a4paper,12pt]{article}
\usepackage{amssymb}
\usepackage{latexsym}
\usepackage{amsmath}
\usepackage{amsthm}
\usepackage{amsfonts}
\usepackage{amssymb}

\usepackage{pstricks}

\newcommand{\rec}[1]{{(\ref{#1})}} 
\newtheorem{Th}{Theorem}[section] 
\newtheorem{Prop}{Proposition}[section]

\newtheorem{Lma}{Lemma}[section]
\newtheorem{Dfi}{Definition}[section] 
\newtheorem{Rm}{Remark}[section] 
 
\newtheorem{Theorem}{Theorem}[section] 
\newtheorem{remark}{Remark}[section]

\newtheorem{Proposition}{Proposition}[section]

\newcommand{\be}{\begin{equation}}
\newcommand{\ee}{\end{equation}}
\newcommand{\bes}{\begin{equation*}}
\newcommand{\ees}{\end{equation*}}

\newcommand{\R}{\mathbb{R}}
\newcommand{\N}{\mathbb{N}}
\newcommand{\C}{\mathbb{C}}
\newcommand{\Z}{\mathbb{Z}}

\def\ti{\tilde}
\def\lf{\left}
\def\rg{\right}

\def\al{\alpha}
\def\la{\lambda}

\def\ep{\varepsilon}

\def\ds{\displaystyle}
\def\ov{\overline}
\def\Om{\Omega}
\def\om{\omega}
\def\p{\partial}

\begin{document}

\title{Horizontal $\alpha$-Harmonic Maps}
\author{ Francesca Da Lio\thanks{Department of Mathematics, ETH Z\"urich, R\"amistrasse 101, 8092 Z\"urich, Switzerland.}  \and   Tristan Riviere$^*$   }
\maketitle
\begin{abstract}
Given a $C^1$ planes distribution $P_T$ on all ${\mathbb R}^m$ we consider {\em horizontal  $\alpha$-harmonic maps}, $\alpha\ge 1/2$,  with respect to such a distribution.
These are maps $u\in H^{\alpha}({\mathbb R}^k,{\R}^m)$ satisfying $P_T\nabla u=\nabla u$ and $P_T(u)(-\Delta)^{\alpha}u=0$  in ${\mathcal D}'({{\mathbb R}}^k).$
  If the distribution of planes is integrable then we recover the classical case of $\alpha$-harmonic maps with values into a manifold. In this paper we shall focus our attention to the case $\alpha=1/2$ in dimension $1$  and $\alpha=2$ in dimension
$2$ and we investigate the regularity of  the {\em horizontal  $\alpha$-harmonic maps}. In both cases we show that such maps satisfy a Schr\"odinger type system with an antisymmetric potential,
that permits us to apply the previous results obtained by the authors in respectively in \cite{Riv} and \cite{DLR2}. Finally we study the regularity of {\em variational $\alpha$-harmonic} maps which are  critical points of $\|(-\Delta)^{\al/2} u\|^2_{L^2}$ 
under the constraint to be tangent (horizontal) to a given planes distribution. We produce a convexification of this variational problem which permits to write it's Euler Lagrange
equations.
\end{abstract}
 \par
\medskip
{\noindent
 {\small {\bf Key words.} Horizontal harmonic map, horizontal fractional harmonic map, sub-riemannian geometry, Schr\"odinger-type PDEs, conservation laws, commutators.}}\par\medskip
{\noindent{\small { \bf  MSC 2010.} 58E20, 35R11, 53C17, 35B65, 35S99, 49Q05.}}
 \tableofcontents 
\medskip
\section{Introduction}
The functions defined on a domain $U\subset {\R}^k$ and which are critical points to the Dirichlet Energy
\[
E(u)=\frac{1}{2}\int_U |\nabla u|^2\ dx^k
\]
where $dx^k$ denotes the Lebesgue measure in ${\R}^k$, satisfy the linear Laplace equation
\[
-\Delta u=0\quad\mbox{ in }{\mathcal D}'(U)
\]
whose solutions are known to be real analytic in any dimension. The result extends of course to maps taking values into
flat euclidian spaces ${\R}^m$. There has been a lot of geometric motivations for studying critical points of the {\it Dirichlet Energy}
among maps forced to take values into a given oriented closed (compact without boundary) sub-manifold $N^n$, that is within the Sobolev
Space defined by
\[
W^{1,2}(U,N^n):=\lf\{ u\in W^{1,2}(U,{\R}^m)\ ;\ u(x)\in N^n\ \mbox{ for a. e. }x\in U\rg\}\quad.
\]
The $C^2$ regularity on $N^n$ is usually assumed in order to ensure at least the Frechet differentiability of $E$ within $W^{1,2}(U,N^n)$. If $\nu$ defines the unit
normal multi-vector to the sub-manifold (under the regularity assumption on $N^n$, we have that $\nu$ is $C^1$), critical points of $E$ satisfy the following Euler Lagrange equation
\be
\label{int-I.1}
\nu(u)\wedge-\Delta u=0\quad \quad\mbox{in }{\mathcal D}'(U)\quad.
\ee
This equation makes sense at the distributional level since the composition of the $C^1$ tensor field $\nu$ with $u$ is in $W^{1,2}$ and $-\Delta u\in H^{-1}(U)$. If $u$ is smooth
equation (\ref{int-I.1}) means that $-\Delta u(x)$  is perpendicular to $T_{u(x)}N^n$ for every $x\in U$ and generalizes the equation of geodesics in $N^n$ for $k=1$ to arbitrary $k$.
Equation (\ref{int-I.1}) is called {\it harmonic map equation} and is often presented in the following equivalent form (see \cite{Riv-cours})
\be
\label{int-I.2}
-\Delta u=\sum_{j=1}^kA(u)(\p_{x_j} u,\p_{x_j} u)
\ee
where $A(z)(X,Y)$ is second fundamental form of $N^n\hookrightarrow {\R}^m$ at the point $z\in N^n$ along the pair of vectors $X,Y\in T_zN^n$.

\medskip

In \cite{DLR1}, the authors initiated the analysis of  {\it $1/2-$harmonic maps} into $N^n$, in connection with the problem of free boundary minimal discs. These maps
are critical points of the fractional energy on ${\R}^k$
\begin{equation}\label{fracenergy}
E^{1/2}(u):=\int_{{\R}^k}|(-\Delta)^{1/4}u|^2\ dx^k
\end{equation}
within
\[
H^{1/2}({\R}^k,N^n):=\lf\{ u\in H^{1/2}({\R}^k,{\R}^m)\ ;\ u(x)\in N^n\ \mbox{ for a. e. }x\in {\R}^k\rg\}\quad.
\]
The corresponding Euler-Lagrange equation is given by
\be
\label{int-I.3}
\nu(u)\wedge(-\Delta)^{1/2} u\quad \quad\mbox{in }{\mathcal D}'({\R}^k)\quad.
\ee
In \cite{DLR1,DLR2}  various regularity results were established for weak solutions to (\ref{int-I.3}) in the critical dimension 1. The proof of these results
where using the existence of special structures in some reformulation of (\ref{int-I.3}) called {\it 3 commutators} which are roughly bilinear pseudo-differential operators satisfying some {\it integrability by compensation}
properties. These results could also be obtained by transforming the a-priori non-local PDE (\ref{int-I.3}) into a local one by performing {\it ad-hoc} extensions and reflections (see \cite{Sche} or \cite{MS}).
More general {\it non-local non-linear elliptic problems} of the form
\be
\label{int-I.4}
\nu(u)\wedge(-\Delta)^{\al} u\quad \quad\mbox{in }{\mathcal D}'({\R}^k)\quad.
\ee
and further generalization have been studied in \cite{Schi1,DL,DLS,}. Observe that, by introducing the field of orthogonal projection $m\times m$ matrices $P_T(z)$ onto the tangent spaces $T_zN^n$
all the above equations can be rewritten in the form
\be
\label{int-I.5}
P_T(u)\,(-\Delta)^\al u=0
\ee
where
\[
\forall\, z\in N^n\quad\forall\, Z\in {\R}^m\quad\quad P_T(z) \,Z:=\sum_{i=1}^mP_T^{ij}(z)\, Z_j\quad\mbox{where }Z=\sum_{j=1}^mZ_j\,\ep_j
\] 
and $(\ep_j)_{j=1\cdots m}$ the canonical basis of ${\R}^m$.

\medskip

The main purpose of the present work is to release the assumption that the field of orthogonal projection $P_T$ is {\it integrable} and associated to a sub-manifold $N^n$
and to consider the equation (\ref{int-I.5}) for a general field of orthogonal projections $P_T$ defined on the whole of ${\R}^m$ and for {\it horizontal maps} $u$ satisfying $P_T(u)\nabla u=\nabla u$.

\medskip

Let $P_T\in C^1({\R}^m,M_m({\R}))$ and $P_N\in C^1({\R}^m,M_m({\R}))$ such that
\be
\label{I.1}\lf\{
\begin{array}{l}
P_T\circ P_T=P_T\quad P_N\circ P_N=P_N\\[3mm]
P_T+P_N=I_m\\[3mm]
\forall\, z\in {\R}^m\quad\forall\, U,V\in T_z{\R}^m\quad <P_TU,P_NV>=0
\end{array}
\rg.
\ee
where $<\cdot,\cdot>$ denotes the standard scalar product in ${\R}^m$. In other words $P_T$ is a $C^1$ map into the orthogonal projections of ${\R}^m$.  For such a distribution of projections $P_T$
we denote by
\[
n:=\mbox{rank}(P_T)\quad.
\]
Such a distribution identifies naturally with the distribution of $n-$planes given by the images of $P_T$ (or the Kernel of $P_T$) and conversely,
any $C^1$ distribution of $n-$dimensional planes defines uniquely $P_T$ satisfying (\ref{I.1}).

For any $\al\ge 1/2$ and for $k\ge 1$
\[
{\frak{H}}^{\al} ({\R}^k):=\lf\{u\in H^{\al}({\R}^k,{\R}^m)\quad;\quad P_N(u)\,\nabla u=0\quad\mbox{ in }{\mathcal D}'({\R}^k)\rg\}
\]
 Observe that this definition makes sense since we have respectively $P_T\circ u\in H^{\al}(S^1,M_m({\R})$ and ${\nabla u}\in H^{\al-1}({\R}^k,{\R}^m)$. We sometimes extend
 this definitions to other domains such as $S^1$ or a general riemannian surface $\Sigma$ for $\al=1$....etc. In the case $\alpha>k/2$ then ${\frak H}^{\al}({\R}^k)$ is a Finsler
 manifold (see Definition 3.8 in \cite{Struwe}).
 \begin{Dfi}
 \label{df-legend-harm}
 Given a $C^1$ plane distribution $P_T$ in ${\R}^m$ satisfying (\ref{I.1}), a map $u$ in the space ${\frak H}^{\al}({\R}^k)$ is called  {\bf {horizontal  $\al$-harmonic} }with respect to $P_T$
 if
 \be
 \label{int-I.6}
 \forall\, i=1\cdots m\quad\quad\sum_{j=1}^mP_T^{ij}(u)(-\Delta)^\al u_j=0\quad\quad\mbox{in }{\mathcal D}'({\R}^k)
 \ee
and we shall use the following notation
\[
 P_T(u)\,(-\Delta)^\al u=0\quad\quad\mbox{in }{\mathcal D}'({\R}^k)\quad.
\] 
\hfill $\Box$
 \end{Dfi}
\noindent{\bf Example}: {\it Horizontal Harmonic Maps} in ${\C}^3$ for the distribution $P_T$  given by
 \be
 \label{int-I.6-a}
P_T(z)\,Z:= Z-|z|^{-2}\ \lf[Z\cdot (z_1,z_2,z_3)\ (z_1,z_2,z_3)+Z\cdot (iz_1,iz_2,i z_3)\ (iz_1,iz_2,i z_3)\rg]
 \ee
are given for instance by the conformal parametrization of {\it Special horizontal Surfaces} in $S^5$ (see \cite{Has}).
  
  \medskip
\begin{Rm}
\label{rm-I.1}
 In \cite{JoY} the authors define a horizontal harmonic map a function which is  in the same time horizontal with respect to the plane distribution \underline{and }
 harmonic.  In \cite{JoY} it is   proved  in particular that in the  case when $P_T$ is issued from a \underbar{riemannian submersion}, also called {\em Carnot-Caratheodory space}, (which is the case of the previous example \rec{int-I.6-a})
 the normal projection of the tension field of any horizontal map in ${\frak H}^1$ is necessary zero. Therefore the horizontal harmonic maps  according to the Definition \ref{df-legend-harm} are both horizontal and harmonic and the two definitions coincide.
  It would be interesting to inquire if such a result holds for $1/2-$harmonic maps.  \hfill $\Box$
\end{Rm}  
  
 \medskip
 
When the plane distribution $P_T$ is {\it integrable} that is to say when
\be
\label{I.2}
\forall \ X,Y\in C^1({\R}^m,{\R}^m)\quad P_N [P_T\, X,P_T\,Y]\equiv 0
\ee
where $[\cdot,\cdot]$ denotes the Lie Bracket of vector-fields, using Fr\"obenius theorem the plane distribution correspond to
the tangent plane plane distribution of a $n-$dimensional {\it foliation } ${\mathcal F}$, (see e.g \cite{Lang}). A smooth map $u$ in ${\frak H}^{\al}({\R}^m)$ takes values  everywhere
into a {\it leaf } of ${\mathcal F}$ that we denote $N^n$ and we are back to the classical theory of harmonic maps into
manifolds. Observe that  our definition includes the case of   $\alpha$-harmonic maps  with  values into a sub-manifold of the euclidean space and horizontal with respect to a planes distribution in this  sub-manifold.
Indeed it is sufficient to add to such a distribution the projection to the sub-manifold and extend the all to a tubular neighborhood of the sub-manifold.

\medskip

In the present work we shall mostly focus our attention to the case $\al=1/2$ in critical dimension 1 and $\al=1$ in critical dimension 2. We establish below that
harmonic maps from ${\R}^m$ into plane distributions satisfy an elliptic Schr\"odinger type system with an antisymmetric potential $\Om\in L^{2}({\R}^k,{\R}^k\otimes so(m))$ of the form
\be
\label{int-I.7}
-\Delta u=\Omega(P_T)\cdot\nabla u.
\ee
Hence, following the analysis in \cite{Riv} we deduce in two dimension  the local existence on a disc $D^2$ of $A(P_T)\in L^\infty\cap W^{1,2}(D^2,Gl_m({\R})$ and $B(P_T)\in W^{1,2}(D^2,M_m({\R})$ such that
\be
\label{int-I.8}
\mbox{div}\lf(A(P_T)\,\nabla u\rg)=\nabla B(P_T)\,\nabla^\perp u
\ee
from which the regularity of $u$ can be deduced using Wente's {\it Integrability by compensation} because of the estimate
\be
\label{int-I.9}
\|\nabla B\,\nabla^\perp u\|_{H^{-1}(D^2)}\le \, C\, \|\nabla B\|_{L^2}\ \|\nabla u\|_{L^2}
\ee
One of the main contribution of our present work is to produce conservation laws corresponding to (\ref{int-I.8}) but for general {\it horizontal $1/2-$harmonic maps} : locally, modulo some smoother terms 
coming from the application of non-local operators on cut-off functions, we construct $A(P_T)\in L^{\infty}\cap H^{1/2}(\R,Gl_m({\R})$ and $B(P_T)\in H^{1/2}(\R,Gl_m({\R})$ such that
\be
\label{int-I.9}
(-\Delta)^{1/4}(A(P_T)\, v)=\mathcal{J}(B(P_T),v) +\mbox{cut-off},
\ee
where $v:=(P_T\,(-\Delta)^{1/4}v,{\mathcal R}\, (P_N(-\Delta)^{1/4}v))$ and ${\mathcal R}$ denotes the Riesz operator and ${\mathcal J}$ is a bilinear pseudo-differential operator satisfying
\be
\label{int-I.10}
\|\mathcal{J}(B,v)\|_{H^{-1/2}({\R})}\le C\, \|(-\Delta)^{1/4} B\|_{L^2({\R})}\, \|v\|_{L^2({\R})}.
\ee
These facts imply the following theorem which is one of the main result of the present work.
\begin{Th}
\label{th-I.1}
Let $P_T$ be a $C^1$ distribution of planes (or projections) satisfying (\ref{I.1}). Any map $u\in {\frak H}^{1/2}({\R})$  satisfying
\be
\label{I.4}
P_T(u)\,(-\Delta)^{-1/2}u=0\quad\mbox{ in }{\mathcal D}'({\R})
\ee 
is in $\cap_{\delta<1}C^{0,\delta}({\R})$\quad. \hfill $\Box$
\end{Th}
Solutions to (\ref{I.4}) are of special geometric interest because of the following proposition extending the well known fact in the integrable case
which has been at the origin of the study of $1/2-$harmonic maps (see \cite{DLR2}).
\begin{Prop}
\label{pr-I.1}
An element in ${\frak H}^{1/2}$ satisfying (\ref{I.4}) has an harmonic extension $\ti{u}$ in $D^2$ which is conformal and hence it is the boundary of a minimal disc
whose exterior normal derivative $\p_r\ti{u}$ is orthogonal to the plane distribution given by $P_T$.\hfill $\Box$
\end{Prop}
\noindent{\bf Example} : We consider the field of projections corresponding to (\ref{int-I.6-a}) but in ${\C}^2\setminus\{0\}$ this time. That is
\be
 \label{int-I.10-a}
P_T(z)\,Z:= Z-|z|^{-2}\ \lf[Z\cdot (z_1,z_2)\ (z_1,z_2)+Z\cdot (iz_1,iz_2)\ (iz_1,iz_2)\rg].
 \ee
Example of $u$ satisfying (\ref{I.4}) is given by solutions to the system
\be
\label{IV.1}
\lf\{
\begin{array}{l}
\ds\frac{\p \ti{u}}{\p r}\in\ \mbox{Span}\lf\{u, i\,u\rg\}\quad\mbox{ a. e.}\\[5mm]
\ds u\cdot\frac{\p u}{\p \theta}=0\quad\mbox{ a. e.}\\[5mm]
\ds i\,u\cdot\frac{\p u}{\p \theta}=0\quad\mbox{ a. e.}
\end{array}
\rg.
\ee
where $\ti{u}$ denotes the harmonic extension of $u$ which happens to be conformal due to proposition~\ref{pr-I.1} and define a minimal disc. An example of such a map is given by
\begin{equation}\label{exhalf}
u(\theta):=\frac{1}{\sqrt{2}}(e^{i\theta}, e^{-i\theta})\quad\mbox{ where }\quad\ti{u}(z,\ov{z})=\frac{1}{\sqrt{2}}(z,\ov{z}).
\end{equation}
Observe the solution in \rec{exhalf} is also an $1/2$-harmonic into $S^3$ and 
it would be interesting to investigate whether this is the unique solution.\footnote{A uniqueness result of that form can be obtained from \cite{FS} in the integrable case when
\[
P_T(z)\,Z:= Z-|z|^{-2}\ Z\cdot (z_1,z_2)\ (z_1,z_2)
\] } to (\ref{I.4}) for $P_T$ given by (\ref{int-I.10-a}) modulo the composition with M\"obius transformations of the form
\[
e^{i\theta}\longrightarrow e^{i\sigma_0}\,\frac{e^{i\theta}-a}{1-\ov{a}\, e^{i\theta}}
\]
where $\sigma_0\in {\R}$, $a\in {\C}$ and $|a|<1$.

Despite the geometric relevance of equations (\ref{int-I.6}) in the non-integrable case, it is however a-priori not the {\it Euler-Lagrange} equation of the variational
problem consisting in finding the critical points of $\|(-\Delta)^{\al/2} u\|^2_{L^2}$ within ${\frak H}^\al$ when $P_T$ is not satisfying (\ref{I.2}). 
This can be seen in the particular case where $\al=1$ where the critical points to the {\it Dirichlet Energy} have been extensively studied in relation with the
computation of {\it normal geodesics} in sub-riemannian geometry.
We then introduce the following definition
\begin{Dfi}
\label{df-I.1}
A map $u$ in ${\frak H}^\al$ is called {\bf variational $\al-$harmonic} into the plane distribution $P_T$  if it is a critical point of the $\|(-\Delta)^{\al/2} u\|^2_{L^2}$
within variations in ${\frak H}^\al$ i.e. for any $u_t\in C^1((-1,1),{\frak H}^\al)$ we have
\[
\lf.\frac{d}{dt}\|(-\Delta)^{\al/2} u_t\|^2_{L^2}\rg|_{t=0}=0.
\] \hfill $\Box$
\end{Dfi}
Example of variational harmonic maps from $S^1$ into plane distribution is given by the sub-riemannian geodesics.\par

The last goal of the present work is to establish an {\it Euler Lagrange} equations characterizing ``variational $\al-$harmonic into the plane distribution $P_T$''
for $\al=1$ and $\al=1/2$ and to study the regularity of these solutions. This is done using a convexification of the variational problem following the spirit
of the approach introduced by Strichartz in \cite{Stri} for {\it normal geodesics} in sub-riemannian geometry. We prove in particular for the case $\al=1/2$ that the smooth critical points
of
\begin{equation}\label{convexlagr}
\begin{array}{l}
\ds{\mathcal L}^{1/2}(u,\xi):=\int_{S^1}\frac{|(-\Delta)^{-1/4}_0(P_T(u)\xi)|^2}{2}\ d\theta\\[5mm]
\ds\quad\quad\quad-\int_{S^1}\lf<(-\Delta)^{-1/4}_0(P_T(u)\xi),(-\Delta)^{-1/4}_0\lf(P_T(u)\frac{du}{d\theta}\rg)\rg>\ d\theta\\[5mm]
\ds \quad\quad\quad-\int_{S^1}\lf<(-\Delta)^{-1/4}_0(P_N(u)\xi),(-\Delta)^{-1/4}_0\lf(P_N(u)\frac{du}{d\theta}\rg)\rg>\ d\theta
\end{array}
\end{equation}
in the co-dimension $m$ \underbar{Hilbert subspace} of $H^{1/2}(S^1,{\R}^m)\times H^{-1/2}(S^1,{\R}^m)$ given by
\[
{\mathfrak E}:=\lf\{
\begin{array}{c}
\ds(u,\xi)\in H^{1/2}(S^1,{\R}^m)\times H^{-1/2}(S^1,{\R}^m)\quad\mbox{s. t. }\\[5mm]
\ds\lf(P_N(u),\frac{du}{d\theta}\rg)_{H^{1/2},H^{-1/2}}=0\\[5mm]
\ds(-\Delta)_0^{-1/4}(P_T(u)\xi)\in L^2(S^1)\quad\mbox{ and }\quad(-\Delta)^{-1/4}_0\lf(P_T(u)\frac{du}{d\theta}\rg)\in L^2(S^1)\quad 
\end{array}\rg\}
\]
at the point where the constraint $\lf(P_N(u),\frac{du}{d\theta}\rg)_{H^{1/2},H^{-1/2}}$ is {\bf  non-degenerate } are ``variational $1/2-$harmonic"  into the plane distribution $P_T$
in the sense of definition~\ref{df-I.1}.  It   remains open the regularity of critical points of \rec{convexlagr} or even of the $1/2$ energy \rec{fracenergy} in ${\frak{H}}^{1/2}$ in the case  when the constraint $\lf(P_N(u),\frac{du}{d\theta}\rg)_{H^{1/2},H^{-1/2}}$ is degenerate.

The paper is organized as follows. In Section 2 we prove the regularity of horizontal Harmonic maps in $2$ dimension. In Section 3 we recall some commutators estimates,
we find conservation laws associated to nonlocal Schr\"odinger type systems with antisymmetric potentials. In Section 4 we deduce  Theorem \ref{th-I.1} from the results obtained in Section 2 and Theorem \ref{pr-I.1}.
In Section 5 we  find the Euler Lagrange equation associated to the Lagrangian \rec{convexlagr} and we show that smooth critical points of \rec{convexlagr}  are actually variational harmonic and  $1/2-$harmonic maps  into a plane distribution $P_T$ .\par
We finally mention that in the case of $1/2-$harmonic maps in dimension $1$ we will consider as domain of definition  indifferently either the real line $\R$ or the circle $S^1.$
 
\medskip

\section{Regularity of  Horizontal Harmonic Maps in $2$-D }

We prove the following theorem.
\begin{Th}
\label{th-reg-harm-map}
Let $P_T$ be a $C^1$ map satisfying (\ref{I.1}). Any map $u\in {\frak H}^1(D^2)$  satisfying
\be
\label{IV-az-1}
P_T(u)\,(-\Delta u)=0\quad\mbox{ in }{\mathcal D}'(D^2)
\ee 
is in $\cap_{\delta<1}C^{0,\delta}_{loc}(D^2)$\quad.
\hfill $\Box$
\end{Th}
\noindent{\bf Proof of theorem~\ref{th-reg-harm-map}.}
We have
\be
\label{IV-az-2}
\begin{array}{l}
-\Delta u=\mbox{div}(\nabla u)=\mbox{div}(P_T(u)\,\nabla u)=P_T(u)\,(-\Delta u)+\nabla(P_T(u))\cdot\nabla u\\[5mm]
\quad=\nabla(P_T(u))\cdot\nabla u=\nabla(P_T(u))\,P_T(u)\cdot\nabla u=-\nabla(P_N(u))\,P_T(u)\cdot\nabla u
\end{array}
\ee
Observe that in one hand
\be
\label{IV-az-3}
\nabla(P_N(u))\,P_T(u)+P_N(u)\,\nabla P_T(u)=0
\ee
and in the other hand
\be
\label{IV-az-4}
(P_N(u)\,\nabla P_T(u))^t=\nabla P_T(u)\,P_N(u)
\ee
Hence combining (\ref{IV-az-2}), (\ref{IV-az-3}) and (\ref{IV-az-4}) together with the fact that $P_N(u)\nabla u\equiv 0$ we obtain
\be
\label{IV-az-5}
-\Delta u=\lf[P_N(u)\,\nabla P_T(u)-(P_N(u)\,\nabla P_T(u))^t\rg]\cdot\nabla u
\ee
Denote $\Om:=\lf[P_N(u)\,\nabla P_T(u)-(P_N(u)\,\nabla P_T(u))^t\rg]$. We have $\Om\in L^2(D^2,so(m)\otimes{\R}^2)$. We can then apply the main result
in \cite{Riv} and deduce theorem~\ref{th-reg-harm-map}.\hfill $\Box$

\section{3-Commutators, Antisymmetry and Conservation Laws in $1$-D}
 
\label{anti}
\subsection{A Regularity Result for Solutions to Linear Pseudo-Differential Equations involving Projections}

Denote  ${\mathcal R}$ is the Riesz operator given by
\[
{\mathcal R}\ :\ f=\sum_{n\in {\Z}}f_n\, e^{i\,n\,\theta}\ \longrightarrow\ {\mathcal R}f:=i\,\sum_{n\in {\Z}^\ast} \mbox{sgn}(n) \,f_n\, e^{i\,n\,\theta}
\]
The following lemma is a straightforward consequence of the classical Coifmann Rochberg and Weiss integrability by compensation (see \cite{CRW}).
\begin{Th}
\label{lm-integ-comp}
Let $m\in {\N}^\ast$, then there exists $\delta>0$ such that for any $P_T,\,P_N\in H^{1/2}(S^1,M_m({\R}))$ satisfying
\be
\label{IV-as.19}\lf\{
\begin{array}{l}
P_T\circ P_T=P_T\quad P_N\circ P_N=P_N\\[3mm]
P_T+P_N=I_m\\[3mm]
\mbox{ for a. e. }\, e^{i\theta}\in S^1\quad\forall\, U,V\in {\R}^m\quad <P_T(\theta)U,P_N(\theta)V>=0
\end{array}
\rg.
\ee
and
\be
\label{IV-as.20}
\int_{S^1}|(-\Delta)^{1/4}P_T|^2\ d\theta<\delta
\ee
then for any $p>1$ and for any $f\in L^p(S^1)$ with $\int_{S^1}f(\theta)\, d\theta=0$
\be
\label{IV-as.21}
\lf(P_T+P_N\,{\mathcal R}\rg)\, f=0\quad\Longrightarrow\quad f=0 
\ee
\end{Th}
We are going to extend the previous theorem to negative Sobolev Spaces.
\begin{Th}
\label{lm-integ-comp-H^{-1/2}}
Let $m\in {\N}^\ast$, then there exists $\delta>0$ such that for any $P_T,\,P_N\in H^{1/2}(S^1,M_m({\R}))$ satisfying
\be
\label{IV.as-19-a}\lf\{
\begin{array}{l}
P_T\circ P_T=P_T\quad P_N\circ P_N=P_N\\[3mm]
P_T+P_N=I_m\\[3mm]
\mbox{ for a. e. }\, e^{i\theta}\in S^1\quad\forall\, U,V\in {\R}^m\quad <P_T(\theta)U,P_N(\theta)V>=0
\end{array}
\rg.
\ee
and
\be
\label{IV.as-20-a}
\int_{S^1}|(-\Delta)^{1/4}P_T|^2\ d\theta<\delta
\ee
then  for any $f\in H^{-1/2}(S^1)$ with $<1,f>_{H^{1/2},H^{-1/2}}=0$
\be
\label{IV.as-21-a}
\lf(P_T+P_N\,{\mathcal R}\rg)\, f=0\quad\Longrightarrow\quad f=0 
\ee
\end{Th}
As we will see in the next subsections, the results on anti-commutators in \cite{CRW} does not apply to the negative Sobolev Spaces
and we are going to make use of {\it integrability by compensation} results for the so called {\it 3-commutators} introduced in \cite{DLR1} combined
with Gauge theoretic arguments exploiting the antisymmetry of some terms in the spirit of \cite{DLR2}. 

\medskip

The uniqueness result~\ref{lm-integ-comp-H^{-1/2}} under small energy assumptions implies the following regularity results
\begin{Th}
\label{th-reg-comp}
Let $m\in {\N}^\ast$ and $P_T,\,P_N\in H^{1/2}(S^1,M_m({\R}))$ satisfying
\be
\label{IV.as-19-aa}\lf\{
\begin{array}{l}
P_T\circ P_T=P_T\quad P_N\circ P_N=P_N\\[3mm]
P_T+P_N=I_m\\[3mm]
\mbox{ for a. e. }\, e^{i\theta}\in S^1\quad\forall\, U,V\in {\R}^m\quad <P_T(\theta)U,P_N(\theta)V>=0
\end{array}
\rg.
\ee
then  for any $f\in H^{-1/2}(S^1)$ with $<1,f>_{H^{1/2},H^{-1/2}}=0$ and satisfying
\be
\label{IV.as-21-a}
\lf(P_T+P_N\,{\mathcal R}\rg)\, f=0
\ee
we have $f\in L^p(S^1)$ for any $p<+\infty$.\hfill$\Box$
\end{Th}
\subsection{Multiplying 3-Commutators}
In this section we recall   regularity properties of some commutators we have introduced in \cite{DLR1,DLR2}, called 3-commutators, and establish
almost stability properties of 3-commutators under multiplication.

\medskip

We introduce the following {\em commutators}:
 
\begin{equation}\label{opT}
T(Q,v):=(- \Delta)^{1/4}(Qv)-Q(- \Delta)^{1/4} v+ (- \Delta)^{1/4} Q v\,\end{equation}
and
\begin{equation}
\label{opS}
S(Q,v):=(-\Delta)^{1/4}[Qv]-{\cal{R}} (Q{\cal{R}}(-\Delta)^{1/4} v)+{\cal{R}}((-\Delta)^{1/4} Q{\cal{R}} v )
\end{equation}
 \begin{equation}\label{opF}
 F(Q,v):={\mathcal{R}}[Q]{\mathcal{R}}[v]-Qv.
 \end{equation}
 \begin{equation}\label{opOm}
 \Lambda(Q,v):=Qv+{\mathcal{R}}[Q{\mathcal{R}}[v]]\,.
 \end{equation}
 In \cite{DLR1,DLR2} the authors obtained the following estimates.
  \begin{Theorem}\label{comm1}
  Let $v\in L^2(\R),$ $Q\in \dot{H}^{1/2}( \R)$ . Then $T(Q,v),S(Q,v)\in {{H}}^{-1/2}( \R)$  and
  \begin{equation}
\label{zz7tris}
\|T(Q,v)\|_{{{H}}^{-1/2}( \R)}\le C\ \|Q\|_{\dot{H}^{1/2}( \R)}\|v\|_{L^{2,\infty}( \R)}\,;
\end{equation}
\begin{equation}
\label{zz8tris}
\|S(Q,v)\|_{{{H}}^{-1/2}( \R)}\le C\ \|Q\|_{\dot{H}^{1/2}( \R)} \|v\|_{L^{2,\infty}( \R)}\,.
\end{equation}
\end{Theorem}
Actually in \cite{DL2} we improve the estimates on the operators $T,S$.
 \begin{Theorem}\label{comm2}
Let Let $v\in L^2(\R),$ $Q\in \dot{H}^{1/2}( \R)$. Then $T(Q,v),S(Q,v)\in {\cal{H}}^{1}( \R)$  and
\begin{equation}\label{commest2}
\|T(Q,v)\|_{{\cal{H}}^{1}( \R)}\le C\|Q\|_{\dot{H}^{1/2}( \R)}\|v\|_{L^2( \R)}\,.\end{equation}
\begin{equation}\label{commest3} 
\|S(Q,v)\|_{{\cal{H}}^{1}( \R)}\le C\|Q\|_{\dot{H}^{1/2}( \R)}\|v\|_{L^2( \R)}\,.~~~\hfill \Box
\end{equation}
\end{Theorem}

 Theorem \ref{comm1} is actually consequence of the following estimates for the dual operators of $T$ and $S$.
 
   \begin{Theorem}\label{commT*}
Let $u,Q\in \dot H^{1/2}( \R)$, denote
$$
T^*(Q,u)=(-\Delta)^{1/4}(Q(-\Delta)^{1/4} u)-(-\Delta)^{1/2}(Qu)+(-\Delta)^{1/4}(((-\Delta)^{1/4} Q) u)\,.$$
then $T^*(Q,u)\in {\cal{H}}^{1}( \R)$ and
\begin{equation}\label{zz9bis}
\|T^*(Q,u)\|_{{\cal{H}}^{1}( \R)}\le C\|Q\|_{\dot{H}^{1/2}( \R)}\|u\|_{\dot{H}^{1/2}( \R)}\,.\hfill\Box\end{equation}
\end{Theorem}
\begin{Theorem}\label{commDLR4bis}
Let $u,Q\in \dot H^{1/2}( \R)$, denote $$
S^*(Q,u)=(-\Delta)^{1/4}(Q(-\Delta)^{1/4} u)-\nabla(Q {\cal{R}}u)+{\cal{R}}(-\Delta)^{1/4}((-\Delta)^{1/4} Q {\cal{R}}u)\,.$$
 Then $S^*(Q,u)\in {\cal{H}}^1( \R)$ and 
\begin{equation}\label{zz10bis}
\|S^*(Q,u)\|_{{\cal{H}}^{1}( \R)}\le C\|Q\|_{\dot{H}^{1/2}( \R)}\|u\|_{\dot{H}^{1/2}( \R)}\,.\hfill \Box\end{equation}
\end{Theorem}
Finally we have 
  \begin{Theorem}
  Let $P,Q\in \dot H^{1/2}( \R)$, denote $$
 \bar{T}(P,Q)= (-\Delta)^{1/4}P{\mathcal{R}}[ (-\Delta)^{1/4} Q]+(-\Delta)^{1/4}[{\mathcal{R}} (-\Delta)^{1/4}[P] Q]-\nabla [PQ].$$
 Then $\bar{T}(P,Q)\in {\cal{H}}^1( \R)$ and 
\begin{equation}\label{estbarT}
\|\bar{T}(P,Q)\|_{{\cal{H}}^{1}( \R)}\le C\|Q\|_{\dot{H}^{1/2}( \R)}\|P\|_{\dot{H}^{1/2}( \R)}\,.\hfill \Box\end{equation}
\end{Theorem}
 \begin{Theorem}\label{estF}
For $f,v\in L^2$ it holds
\begin{equation}\label{estFinftybis}
\|F(f,v)\|_{  H^{-1/2}(\R)}\le C \|f\|_{L^2(\R)}\|v\|_{L^{2,\infty}(\R)}\,.
\end{equation}
and 
\begin{equation}\label{estFhardy}
\|F(f,v)\|_{  {\cal{H}}^1(\R)}\le C \|f\|_{L^2(\R)}\|v\|_{L^{2}(\R)}\,.
\end{equation}
\end{Theorem}
 
\begin{Theorem}\label{estOm}
For  $Q\in \dot H^{1/2}( \R)$, $v\in L^2( \R)$  it holds
\begin{equation}\label{estFhardy}
\|\Lambda(Q,v)\|_{ L^{2,1}( \R)}\le C \|Q\|_{H^{1/2}( \R)}\|v\|_{L^{2}( \R)}\,.
\end{equation}
\end{Theorem}
Actually the estimate \rec{estFhardy} is a consequence of the Coifman-Rochberg- Weiss estimate \cite{CRW}.

 Next we  prove a sort of stability of  of the operators $T,S,F$ with respect to the multiplication  by a  function $P\in H^{1/2}(\R)\cap L^{\infty}(\R).$  Roughly speaking 
 if we multiply them by  a  function $P\in H^{1/2}(\R)\cap L^{\infty}(\R)$ we get a decomposition into the sum of a function in the Hardy Space and a term which is the product of function in $L^{2,1}$ by one in $L^2.$
 
 \begin{Theorem}{ {\bf [Multiplication of $F$ by $P\in H^{1/2}(\R)\cap L^{\infty}(\R)$ ]}}\label{multF}
 Let  $P\in  H^{1/2}(\R)\cap L^{\infty}(\R)$ and $f,v\in L^2(\R)$. Then
 \begin{equation}
PF(f,v)= \underbrace{F(P{\mathcal{R}}[f],{\mathcal{R}} [v])}_{\in {\cal{H}}^1(\R)}-\underbrace{\Lambda(P,f)}_{\in L^{2,1}}\,v
\end{equation}
 \end{Theorem}
 {\bf Proof of Theorem \ref{multF}.}
 We have
 \begin{eqnarray*}
 PF(f,v)&=&P{\mathcal{R}}[f]\,{\mathcal{R}}[v]-Pfv\\
 &=& P{\mathcal{R}}[f]\,{\mathcal{R}}[v]+{\mathcal{R}}[P{\mathcal{R}}[f]] v
 -{\mathcal{R}}[P{\mathcal{R}}[f]] v- Pfv\\
 &=&F(P{\mathcal{R}}[f],{\mathcal{R}} [v])- {\Lambda(P,f)}v.
 \end{eqnarray*}
  The conclusion follows from Theorem \ref{estF} and Theorem \ref{estOm}.\hfill$\Box$

\begin{Theorem}{ {\bf [Multiplication of $T$ by $P\in H^{1/2}(\R)\cap L^{\infty}(\R)$]}}\label{multT}

Let  $P,Q\in  H^{1/2}(\R)\cap L^{\infty}(\R)$ and $v\in L^2(\R)$. Then
\begin{equation}
PT(Q,v)=J_{T}(P,Q,v)+{{\mathcal{A}}_{T}(P,Q)}  v,
\end{equation}
where
$$
 {\mathcal{A}}_T(P,Q)=(-\Delta)^{1/4}T^*(P,Q)=P(-\Delta)^{1/4}[Q]+ (-\Delta)^{1/4}[P] Q-(-\Delta)^{1/4}[PQ]\in L^{2,1}$$
 with
 \begin{equation}\label{calA}
 \|{\mathcal{A}}_T(P,Q)\|_{L^{2,1}}\leq C\|(-\Delta)^{1/4}[P] \|_{L^2}\|(-\Delta)^{1/4}[Q] \|_{L^2},\end{equation}
 and
$$
J_{T}(P,Q,v):= T(PQ,v)-T(P,Qv)\in  {\cal{H}}^1(\R)$$ with
\begin{equation}\label{J}
\|J_{T}(P,Q,v)\|_{{\cal{H}}^1(\R)}\le C \left(\|(-\Delta)^{1/4}[P] \|_{L^2}+\|(-\Delta)^{1/4}[Q] \|_{L^2 }\right)\|v\|_{L^2}.\end{equation}
\end{Theorem}
{\bf Proof of Theorem \ref{multT}.}
We have
 \begin{eqnarray*}
 PT(Q,v)&=&P(- \Delta)^{1/4}[Qv]-PQ(- \Delta)^{1/4} [v]+ P(- \Delta)^{1/4} [Q] v\\
 &=&\{P (- \Delta)^{1/4} [Q]-(- \Delta)^{1/4} [PQ]+(- \Delta)^{1/4}[P]Q\}v\\
 &+&(- \Delta)^{1/4} [PQv]-PQ(- \Delta)^{1/4} v+(- \Delta)^{1/4} [PQ]v\\&-&\left((- \Delta)^{1/4} [PQv]+P(- \Delta)^{1/4}(Qv)-(- \Delta)^{1/4}[P]Qv\right)
 \\
 &=&(-\Delta)^{-1/4}[T^*(P,Q)] v+T(PQ,v)-T(P,Qv).
 \end{eqnarray*}
 Finally  the estimates \rec{calA}, \rec{J} follow from Theorem \ref{commT*} and Theorem \ref{comm2}.\hfill$\Box$
 
 \medskip
 
Now we consider the operator $S$. We first observe that given $Q\in H^{1/2}$ and $v\in L^2$ we have the following decomposition

\begin{eqnarray}\label{decS}
{\mathcal{R}}[S(Q,v)]&=&\tilde S(Q,v)-{\mathcal{R}} (-\Delta)^{1/4} [Q] v-(-\Delta)^{1/4} Q{\mathcal{R}}v\\[5mm]
&=& \tilde S(Q,v)+F({\mathcal{R}} (-\Delta)^{1/4} [Q],v).\nonumber
\end{eqnarray}

where
 $$
\tilde S(Q,v)={\mathcal{R}}(-\Delta)^{1/4}[Qv]+  Q{\cal{R}} (-\Delta)^{1/4}[v]+ {\mathcal{R}} (-\Delta)^{1/4} [Q] v.
$$
From Theorems \ref{comm2} and \ref{estF} it follows that 
 We observe that $\tilde S(Q,v)\in {\cal{H}}^{1}$ and
$$\|\tilde S(Q,v)\|_{{\cal{H}}^{1}}\le \|v\|_{L^2}\|Q\|_{H^{1/2}}.$$
 
 \begin{Theorem}{{\bf [Multiplication of ${\cal{R}}S$ by a rotation $P\in H^{1/2}(\R)\cap L^{\infty}(\R)$]}}\label{multS}
Let  $P,Q\in  H^{1/2}(\R)\cap L^{\infty}(\R)$ and $v\in L^2(\R)$. Then
\begin{eqnarray}\label{decS}
P{\mathcal{R}}[S(Q,v)]&=&  {\mathcal{A}}_S(P,Q) v+J_{S}(P,Q,v)\end{eqnarray}
where
$$ {\mathcal{A}}_S(P,Q):=(-\Delta)^{-1/4}[\bar{T}(P,Q)]+ \Lambda(P, {\mathcal{R}} (-\Delta)^{1/4} [Q])\in L^{2,1}.$$
with
 $$
 \| {\mathcal{A}}_S(P,Q)\|_{L^{2,1}}\leq C\|(-\Delta)^{1/4}[P] \|_{L^2}\|(-\Delta)^{1/4}[Q] \|_{L^2},$$
 and
 $$
J_{S}(P,Q,v):= \tilde S(PQ,v)-\tilde S(P,Qv)+F({\mathcal{R}}[P(-\Delta)^{1/4} [Q]],v)\in  {\cal{H}}^1(\R)$$ with
$$
\|J_{S}(P,Q,v)\|_{{\cal{H}}^1(\R)}\le C \left(\|(-\Delta)^{1/4}[P] \|_{L^2}+\|(-\Delta)^{1/4}[Q] \|_{L^2 }\right)\|v\|_{L^2}.$$

  \end{Theorem}
 \par
  {\bf Sketch of Proof.}
Let  $P\in  H^{1/2}(\R)\cap L^{\infty}(\R)$, then
\begin{eqnarray}
P{\mathcal{R}}[S(Q,v)]&=&\tilde S(PQ,v)-\tilde S(P,Qv)\\[5mm]
&+&\{P{\mathcal{R}} (-\Delta)^{1/4} Q+{\mathcal{R}} (-\Delta)^{1/4}[P] Q-{\mathcal{R}} (-\Delta)^{1/4}[PQ]\}\, v\nonumber\\[5mm] & -&
P[{\mathcal{R}} (-\Delta)^{1/4} [Q] v-(-\Delta)^{1/4} Q{\mathcal{R}}v].\nonumber
\end{eqnarray}
Next we estimate  the term $P[{\mathcal{R}} (-\Delta)^{1/4} [Q] v-(-\Delta)^{1/4} Q{\mathcal{R}}v]$
\begin{eqnarray}
 P[{\mathcal{R}} (-\Delta)^{1/4} [Q] v+(-\Delta)^{1/4} Q{\mathcal{R}}v]&=&\underbrace{\{P{\mathcal{R}} (-\Delta)^{1/4} [Q] -{\mathcal{R}} [P (-\Delta)^{1/4} [Q] ]]\}}_{\in L^{2,1}}v\nonumber\\[5mm]
 &+&\underbrace{{{\mathcal{R}} [P (-\Delta)^{1/4} [Q] ]\,v+P(-\Delta)^{1/4} Q\,{\mathcal{R}}[v]}}_{\in {{\cal{H}}^{1}}}.
 \end{eqnarray}
 Therefore we can write
 \begin{eqnarray}\label{decS}
P{\mathcal{R}}[S(Q,v)]&=&  {\mathcal{A}}_S(P,Q) v\\[5mm] &+&\underbrace{\tilde S(PQ,v)-\tilde S(P,Qv)+F({\mathcal{R}}[P(-\Delta)^{1/4} [Q]],v)}_{\in{\cal{H}}^1}\nonumber
\end{eqnarray}
where
$$ {\mathcal{A}}_S(P,Q):=(-\Delta)^{-1/4}[\bar{T}(P,Q)]+ \Lambda(P, {\mathcal{R}} (-\Delta)^{1/4} [Q]).~~~\hfill\Box$$

 \begin{remark}{\rm 
  We   mention without entering into the details that in $2$-D the Jacobian $J(a,b)=\nabla(a)\,\nabla^{\perp}(b)$ satisfies  a   stability property enjoyed by the operators \rec{opT}, \rec{opS}, \rec{opF} with respect to the multiplication 
by   $P\in W^{1,2}(\R^2)\cap L^{\infty}(\R^2)$ as well. More precisely we may define the following two zero order pseudo-differential operators:
$\mbox{Grad}( X):=\nabla \mbox{div}\Delta^{-1}(X)$, $\mbox{Rot} (Y)=\nabla^{\perp} \mbox{curl}\Delta^{-1}(Y)$. If $a,b\in W^{1,2}(\R^2)$ and  $P\in W^{1,2}(\R^2)\cap L^{\infty}(\R^2)$ then

\begin{eqnarray}\label{jac}
J(a,b)&=&\nabla(a)\,\nabla^{\perp}(b)\\&=&\mbox{Grad}(\nabla(a))\,\mbox{Rot}(\nabla^{\perp}(b))-\mbox{Rot}(\nabla(a))\,\mbox{Grad}(\nabla^{\perp}(b));\nonumber
\end{eqnarray}
and 
\begin{eqnarray}\label{jac}
P\,J(a,b)&=&P\,\nabla(a)\,\nabla^{\perp}(b)\\&=&
\underbrace{[P\mbox{Grad}(\nabla(a)) -\mbox{Grad}(P\nabla(a))]}_{\in L^{2,1}(\R^2)}\,\mbox{Rot}(\nabla^{\perp}(b))\nonumber\\&+&
\underbrace{\mbox{Grad}(P\nabla(a))\,\mbox{Rot}(\nabla^{\perp}(b))-\mbox{Rot}(P\nabla(a))\,\mbox{Grad}(\nabla^{\perp}(b))}_{\in{\cal{H}}^1(\R^2)}.\nonumber
\end{eqnarray}}

\end{remark}

\subsection{Conservation Laws for Fractional Schr\"odinger type PDEs with Antisymmetric Potentials.}

The aim of this part is to construct conservation laws for fractional Schr\"odinger type PDEs with antisymmetric potentials. More precisely we are going to consider a nonlocal system of the form
\begin{equation}\label{modelsystem}
(-\Delta)^{1/4} v=\Omega_0 v+\Omega_1 v+{\cal{Z}}(Q,v)+g(x)\end{equation}
where $v\in L^2(\R)$, $Q\in  H^{1/2}(\R)$, ${\cal{Z}}\colon  H^{1/2}(\R)\times L^2(\R)\to {\cal{H}}^1(\R)$ is 
a linear combination of the operators \rec{opF}, \rec{opT} and \rec{opS}  introduced in the previous section, $\Omega_0\in L^2(\R,so(m))$, $\Omega_1\in L^{2,1}(\R)$, $g(x)$ is a tempered distribution.\par
\begin{Theorem}\label{conservation}
 Let $v\in L^{2}(\R,\R^m)$ be a solution of \rec{modelsystem},
 where $\Omega_0\in L^2(\R,so(m))$, $\Omega_1\in L^{2,1}(\R)$,  ${\cal{Z}}$ is a linear combination of the operators  \rec{opF}, \rec{opT} and \rec{opS}, ${\cal{Z}}(Q,v)\in {\cal{H}}^1$ for every $Q\in{H^{1/2}},$ $v\in L^2$ with
\begin{eqnarray*}
\|{\cal{Z}}(Q,v))\|_{{\cal{H}}^1}&\le& C \| Q\|_{{H^{1/2}}}\|v\|_{L^2}
\end{eqnarray*}
There exists $\varepsilon_0>0$ such that if $\|\Omega_0\|_{L^2}<\varepsilon_0$, then 
 there exist    $A=A(\Omega_0,\Omega_1,Q)\in H^{1/2}(\R,M_m({\R})))$ and an operator $ B=B(\Omega_0,\Omega_1Q)\in  H^{1/2}(\R)$ such that
\begin{eqnarray}
\|A\|_{ H^{1/2}}+\|B\|_{ H^{1/2}}&\le& C(\|\Omega_0\|_{L^2}+\|\Omega_0\|_{L^{2,1}} +\| Q\|_{{H^{1/2}}})\\
dist(A, SO(m))&\le& C(\|\Omega_0\|_{L^2}+\|\Omega_0\|_{L^{2,1}} +\| Q\|_{{H^{1/2}}})
\end{eqnarray}
and 
\begin{equation}\label{conslaw}
(-\Delta)^{1/4}[Av]={\cal{J}}(B,v)+Ag,
\end{equation}
where ${\cal{J}}$ is a linear operator in $B,v$, ${\cal{J}}(B,v)\in {\cal{H}}^1(\R)$ and 
\begin{equation}\label{J}
\|{\cal{J}}(B,v)\|_{{\cal{H}}^1(\R)}\le C \| B \|_{H^{1/2}}\|v\|_{L^2}\,.
\end{equation}
 \end{Theorem}

 {\bf Proof of Theorem \ref{conservation}.} 
 We first observe that since 
 the operator ${\cal{Z}}(Q,v))$   is a linear combination of the operators $F$, $S$ and $T$, it satisfies the following stability  property:    if $Q,P\in{ \dot {H}}^{1/2}(\R)\cap L^{\infty}(\R)$, $v\in L^{2}$ then  
\begin{equation}\label{decH}
P{\cal{Z}}(Q,v)={{{\cal{A}}_{{\cal{Z}}}}}(P,Q)v+J_{{\cal{Z}}}(P,Q,v),
\end{equation}
where
$$
 \|{\cal{A}}_{{\cal{Z}}}(P,Q)\|_{L^{2,1}}\leq C\|(-\Delta)^{1/4}[P] \|_{L^2}\|(-\Delta)^{1/4}[Q] \|_{L^2},$$
 and
 $$
\|{{J}}_{{\cal{Z}}}(P,Q,v)\|_{{\cal{H}}^1(\R)}\le C \left(\|(-\Delta)^{1/4}[P] \|_{L^2}+\|(-\Delta)^{1/4}[Q] \|_{L^2 }\right)\|v\|_{L^2}.$$

{\bf Step 1:} From Theorem 1.2 in \cite{DLR2}  there exists $\varepsilon_0>0$ and $C>0$ such that  
if $\|\Omega_0\|_{L^2}<\varepsilon_0$, then there exists $P=P(\Omega_0)\in{ \dot {H}}^{1/2}(\R, SO(m))$  such that
\begin{equation}\label{cond}\left\{\begin{array}{ll}
(i)~~& \ds P^{-1}(-\Delta)^{1/4} P  -(-\Delta)^{1/4} P^{-1} P=2\,\Omega_0\,;
\\[5mm]
(ii)~~& \|(-\Delta)^{1/4} P\|_{L^2}\le  C\|\Omega_0\|_{L^2} \,.
\end{array}\right.
\end{equation}
Moreover 
\begin{eqnarray*}
P\Omega_0P^{-1}-P^{-1}(-\Delta)^{1/4} P&=&-\frac{\left( P(-\Delta)^{1/4} [P^{-1}]-(-\Delta)^{1/4} [P]P^{-1}\right)}{2}\\
&=&-(-\Delta)^{-1/4}(T^*(P^{-1},P))\in L^{2,1}.\end{eqnarray*}
{\bf Step 2:} {\bf Estimate of $(-\Delta)^{1/4} [Pv]$.}\par
\begin{eqnarray}\label{estdeltaqv}
(-\Delta)^{1/4} [Pv]&=&(-\Delta)^{1/4} [Pv]-P(-\Delta)^{1/4} [v]+(-\Delta)^{1/4} [P]  v\\
&+& P(-\Delta)^{1/4} [v]-(-\Delta)^{1/4} [P]  v\nonumber\\
&=& T(P,v)+P\{\Omega_0 v+\Omega_1 v+{\cal{Z}}(P,v)+g(x)\}-(-\Delta)^{1/4} [P]  v\nonumber
\\
&=&T(P,v)+J_{{\cal{Z}}}(P,Q,v)+{\cal{A}}_{{\cal{Z}}}(P,Q)v+P\Omega_1P^{-1} (P v)\nonumber
\\
&+&[P\Omega_0P^{-1}-(-\Delta)^{1/4} [P]P^{-1}] (P v)+Pg \nonumber\\
&=&\varpi(\Omega_0,\Omega_1,Q)(P v)+J_{T,{\cal{Z}}}(P,v)+Pg \nonumber
\end{eqnarray}
where 
$$\varpi(\Omega_0,\Omega_1,Q)=P\Omega_1P^{-1}+{\cal{A}}_{{\cal{Z}}}(P,Q)+[P\Omega_0P^{-1}-(-\Delta)^{1/4} [P]P^{-1}]\in L^{2,1}$$
with
$$ \|\varpi \|_{L^{2,1}}\leq C (\|(-\Delta)^{1/4}[Q] \|_{L^2}+\|\Omega_0\|_{L^2}+\|\Omega_1\|_{L^{2,1}}),$$
and $J_{T,{\cal{Z}}}(P,Q,v)=T(Q,v)+J_{{\cal{Z}}}(P,Q,v)\in {\cal{H}}^1(\R)$ with
 $$
\|J_{T,{\cal{Z}}}(P,Q,v)\|_{{\cal{H}}^1(\R)}\le C \left(\|(-\Delta)^{1/4}[P] \|_{L^2}+\|(-\Delta)^{1/4}[Q] \|_{L^2 }\right)\|v\|_{L^2}.$$
Moreover the operator ${\cal{J}}$  is linear   and it has   the following property:
 if $ M,P,Q\in{ \dot {H}}^{1/2}(\R)\cap L^{\infty}(\R)$, $v\in L^{2}$ then it holds the following decomposition:
\begin{equation}\label{decH}
MJ_{T,{\cal{Z}}}(P,Q,v)=\omega(M,P,Q,)v+{\cal{G}}(M,P,Q,v),
\end{equation}
with $\omega(M,P,Q,v)\in L^{2,1}(\R)$ and ${\cal{G}}(M,P,Q,v)\in {\cal{G}}^1(\R)$. This decomposition follows from the fact that $J_{T,{\cal{Z}}}(P,Q,v)$ is a linear combination of the operators
$F$, $T$ and $S$.\par
 {\bf Step 3:} Given  ${\mathcal{E}}\in W^{1/2,(2,1)}\cap L^{\infty}$ (that we will chose later in a suitable way), from the above computations it follows that
\begin{eqnarray}\label{estepsv}
(-\Delta)^{1/4} [( Id +{\cal{E}})Pv]&=&T(( Id +{\cal{E}}),Pv)\nonumber\\[5mm]
&+& ( Id +{\cal{E}})(-\Delta)^{1/4} [Pv]-(-\Delta)^{1/4} [( Id +{\cal{E}})]  Pv\nonumber\\[5mm]
&=& T(( Id +{\cal{E}}),Pv)+( Id +{\cal{E}})\{\varpi(\Omega_0,\Omega_1,Q) (Q v)+J_{T,{\cal{Z}}}(P,Q,v)\}\nonumber\\[5mm]
&-&(-\Delta)^{1/4} [{\cal{E}}]  Pv \\[5mm]
&=&T(( Id +{\cal{E}}),Pv)+{\cal{G}}(( Id +{\cal{E}}),P,Q,v)+( Id +{\cal{E}})Pg  \nonumber\\[5mm]
&+&[( Id +{\cal{E}})\varpi(\Omega_0,\Omega_1,Q)+\omega( Id +{\cal{E}}),P,Q)P^{-1}]Pv-(-\Delta)^{1/4} [{\cal{E}}]  Pv.\nonumber
\end{eqnarray}
Set
$\tilde\omega(( Id +{\cal{E}}),\Omega_0,\Omega_1,Q):=[( Id +{\cal{E}})\varpi(\Omega_0,\Omega_1,Q)+\omega (( Id +{\cal{E}}),P,Q)P^{-1}]$. We have $\tilde\omega$ is linear with respect  to ${\cal{E}}$, 
$\tilde\omega(( Id +{\cal{E}}),\Omega_0,\Omega_1,Q)\in L^{2,1}$ and
$$
\|\tilde\omega(( Id +{\cal{E}}),\Omega_0,\Omega_1,Q)\|_{L^{2,1}}\le C\|{\cal{E}}\|_{L^{\infty}} (\|(-\Delta)^{1/4}[Q] \|_{L^2}+\|\Omega_0\|_{L^2}+\|\Omega_1\|_{L^{2,1}}) \,.$$
We choose  ${\cal{E}}$ to  be a solution     in $W^{1/2,(2,1)}\cap L^{\infty}$
   \begin{equation}\label{Eps}
(-\Delta)^{1/4} [{\cal{E}}] =\tilde\omega(Id +{\cal{E}},\Omega_0,\Omega_1,Q).
\end{equation}
Such a solution satisfies
$$
\|{\cal{E}}\|_{L^\infty}\le C\|(-\Delta)^{1/4} [{\cal{E}}] \|_{L^{2,1}}\le C (\|(-\Delta)^{1/4}[Q] \|_{L^2}+\|\Omega_0\|_{L^2}+\|\Omega_1\|_{L^{2,1}}).
$$

By combining  \rec{estepsv} and \rec{Eps} it follows
\begin{equation}\label{estepsvbis}
(-\Delta)^{1/4} [( Id +{\cal{E}})Pv]={\cal{J}}(( Id +{\cal{E}}),P,Q,v)+( Id +{\cal{E}})Pg,
\end{equation}
where 
 $${\cal{J}}(( Id +{\cal{E}}),\Omega_0,\Omega_1,Q,v)=T(( Id +{\cal{E}}),Pv)+{\cal{G}}(( Id +{\cal{E}}),P,Q,v)$$
${\cal{J}}(( Id +{\cal{E}}),\Omega_0,\Omega_1, Q,v)\in  {\cal{H}}^1(\R)$ with
 $$
\|{\cal{J}}(( Id +{\cal{E}}),\Omega_0,\Omega_1, Q,v)\|_{{\cal{H}}^1(\R)}\le C(\|(-\Delta)^{1/4}[Q] \|_{L^2}+\|\Omega_0\|_{L^2}+\|\Omega_1\|_{L^{2,1}}) \|v\|_{L^2}.$$
 We set $A=A(\Omega_0,\Omega_1,Q)=( Id +{\cal{E}})P$ and $B=B(\Omega_0,\Omega_1,Q)=(( Id +{\cal{E}}),P,Q),$ where $P$ satisfies \rec{cond} and ${\cal{E}}$ is a solution of \rec{Eps}. It is evident that
 $( Id +{\cal{E}})$ dependson $\Omega_0,\Omega_1,Q.$
  
 We get
\begin{equation}\label{consfinal}
(-\Delta)^{1/4} [Av]={\cal{J}}(B,v) +Ag
\end{equation}
 We observe that by construction
 we have
 \begin{eqnarray*}
\|A\|_{ H^{1/2}}+\|B\|_{ H^{1/2}}&\le& C(\|\Omega_0\|_{L^2}+\|\Omega_0\|_{L^{2,1}} +\| Q\|_{{H^{1/2}}})\\
dist(A, SO(m))&\le& \|A-P\|_{L^\infty}\le C\|{\cal{E}}\|_{L^\infty}\le C(\|\Omega_0\|_{L^2}+\|\Omega_0\|_{L^{2,1}} +\| Q\|_{{H^{1/2}}})\\
\|{\cal{J}}(B,v)\|_{{\cal{H}}^1(\R)}&\le& C\|B\|_{ H^{1/2}}\|v\|_{L^2}.
\end{eqnarray*}
We can conclude the proof.~~~$\Box$
\begin{Th}
\label{lm-integ-comp-H^{-1/2}bis}
Let $m\in {\N}^\ast$, then there exists $\delta>0$ such that for any $P_T,\,P_N\in H^{1/2}(\R,M_m({\R}))$ satisfying
\be
\label{IV.as-19-abis}\lf\{
\begin{array}{l}
P_T\circ P_T=P_T\quad P_N\circ P_N=P_N\\[3mm]
P_T+P_N=I_m\\[3mm]
\mbox{ for a. e. }\, x\in\R\quad\forall\, U,V\in {\R}^m\quad <P_T(x)U,P_N(x)V>=0
\end{array}
\rg.
\ee
and
\be
\label{IV.as-20-a}
\int_{\R}|(-\Delta)^{1/4}P_T|^2\ d\theta<\delta
\ee
then  for any $f\in H^{-1/2}(\R)$  
\be
\label{IV.as-21-abis}
\lf(P_T+P_N\,{\mathcal R}\rg)\, f=0\quad\Longrightarrow\quad f=0 
\ee
\end{Th}
 {\bf Proof of Theorem \ref{lm-integ-comp-H^{-1/2}bis}.}
 
 We first set $f:=(-\Delta)^{1/2}u$. From \rec{IV.as-21-abis} it follows that
 \begin{equation}\label{harmequ}
\left\{\begin{array}{c}
P_T(-\Delta)^{1/2}u=0\\[5mm]
P_N {\mathcal{R}} (-\Delta)^{1/2}u=0
\end{array}\right.
\end{equation}
The set  $v=(P_T(-\Delta)^{1/4}u,P_N {\mathcal{R}} (-\Delta)^{1/4}u)^t$. In \cite{DLR2} it has been proved that $v$ satisfies a nonlocal Schr\"odinger type system of the form
\rec{modelsystem}
with $g\equiv 0$ 
  $\Omega_0=\Omega_0(P_T)\in L^2(\R,so(\R^m))$  $\Omega_1=\Omega_1(P_T)\in L^{2,1}$,  ${\cal{Z}}(P_T,v)$ is a linear operator in $P_T,v$, ${\cal{Z}}(P_T,v)\in {\cal{H}}^1$ with
\begin{eqnarray*}
\|\Omega_0\|_{L^2}&=&\|\Omega_0(P_T)\|_{L^2}\le C \| P_T\|_{H^{1/2}}\\
\|\Omega_1\|_{L^{2,1}}&=&\|\Omega_1(P_T)\|_{L^{2,1}}\le C \| P_T\|_{{H^{1/2}}}\\
\|{\cal{Z}}(P_T,v))\|_{{\cal{H}}^1}&\le& C \| P_T\|_{{H^{1/2}}}\|v\|_{L^2}
\end{eqnarray*}
(see appendix \ref{rewrite}).
 From Theorem \ref{conservation} it follows that if $\delta$ is small enough then  there exist
  $A=A(P_T) $ and $B=B(P_T)$ such that 
  \begin{equation}\label{consfinalhom}
(-\Delta)^{1/4} [Av]={\cal{J}}(B,v)  
\end{equation}
and 
 \begin{eqnarray}
\|A\|_{ H^{1/2}}+\|B\|_{ H^{1/2}}&\le& C\|P_T\|_{H^{1/2}}\nonumber\\
dist(A, SO(m))&\le&  \le C \|P_T\|_{H^{1/2}}\label{estJ}\\
\|{\cal{J}}(B,v)\|_{{\cal{H}}^1(\R)}&\le& C\|B\|_{ H^{1/2}}\|v\|_{L^2}.\nonumber
\end{eqnarray}
From \rec{consfinalhom} and \rec{estJ} it follows that
\begin{eqnarray}\label{estvfinal}
\|v\|_{L^{2}}&=&\|A^{-1}Av\|_{L^2}\le C\|A^{-1}\|_{L^{\infty}}\|Av\|_{L^{2}}\\
&\le&C \|(-\Delta)^{-1/4}{\cal{J}}(B,v)\|_{L^{2,1}}\le   C\|B\|_{ H^{1/2}}\|v\|_{L^2}\nonumber\\
&\le& C  \|P_T\|_{H^{1/2}}\|v\|_{L^2}\le C  \varepsilon_0 \|v\|_{L^2}. \nonumber
\end{eqnarray} 
Again if $\delta$ is small enough then \rec{estvfinal} yields $v\equiv 0$.
\par
{\bf Claim: } $v\equiv 0~~\Rightarrow f=0.$
\par
{\bf Proof of the Claim.}\par
If $\delta$ is small enough  then $v\equiv 0$ implies that $P_T(-\Delta)^{1/4}u=0$ and $P_N {\mathcal{R}} (-\Delta)^{1/4}u=0.$
Now observe that
\begin{equation}
{\mathcal{R}}[P_N (-\Delta)^{1/4}u]=\underbrace{{\mathcal{R}}[P_N (-\Delta)^{1/4}u]-\underbrace{P_N {\mathcal{R}} (-\Delta)^{1/4}u}_{=0}}_{\in L^{2,1}}\end{equation}
Therefore from Theorem \ref{estF}  and the fact that $P_T(-\Delta)^{1/4}u=0$ it follows that
\begin{eqnarray}\label{claimv}
\| (-\Delta)^{1/4}u\|_{L^2}&=&\| P_N (-\Delta)^{1/4}u]|_{L^2}=\|{\mathcal{R}}[P_N (-\Delta)^{1/4}u]\|_{L^2}\\&\le& C \|P_T\|_{H^{1/2}}\|v\|_{L^2}\le C  \delta \|(-\Delta)^{1/4}u\|_{L^2}.\nonumber
\end{eqnarray} 
  Therefore if $ C  \delta<1$ then \rec{claimv} implies $-\Delta)^{1/4}u\equiv0$ and therefore $f=0$.
  \par
  We can conclude the proof of the claim and of the Theorem \ref{lm-integ-comp-H^{-1/2}bis}.~~$\Box$
  \section{Regularity of Horizontal $1/2-$Harmonic Maps in $1$-D}

\subsection{Proof of Theorem~\ref{th-I.1}.}
 The proof of Theorem \ref{th-I.1} follows by combining Theorem \ref{lm-integ-comp-H^{-1/2}bis} and
  localization arguments used in \cite{DLR2}.~~~$\Box$

\subsection{Proof of Proposition~\ref{pr-I.1}.}

A first proof when $P_T\in C^2({\R}^m)$. In that case we have that $u\in C^{1,\al}(S^1)$. Denote $\ti{u}$ the harmonic extension
of $u$. It is well known that the {\it Hopf differential} of $\ti{u}$
\[
|\p_{x_1}\ti{u}|^2-|\p_{x_2}\ti{u}|^2-\,2\, i\, \lf<\p_{x_1}\ti{u},\p_{x_2}\ti{u}\rg>=f(z)
\]
is holomorphic. Considering on $S^1=\p D^2$
\[
2\,\lf<\p_r\ti{u},\p_\theta\ti{u}\rg>=-\,\sin2\theta\, \lf(|\p_{x_1}\ti{u}|^2-|\p_{x_2}\ti{u}|^2\rg)-\,\cos 2\,\theta\ (-\, 2\, \lf<\p_{x_1}\ti{u},\p_{x_2}\ti{u}\rg>)=-\Im\lf(z^2\,f(z)\rg)
\]
Since $0=P_T(u)\,(-\Delta)^{1/2}u= P_T(u)\,\p_{r}\ti{u}$ and $0=P_N(u)\,\p_{\theta}u=P_N(u)\,\p_{\theta}\ti{u}$ on $\p D^2$ we have that
\[
\Im\lf(z^2\,f(z)\rg)=0\quad \mbox{ on }\p D^2
\]
hence the holomorphic function $z^2\, f(z)$ is equal to a real constant. Since $f(z)$ cannot have a pole at the origin we have that $z^2 f(z)$
is identically equal to zero and hence $\ti{u}$ is conformal.\hfill $\Box$

\section{Variational Harmonic Maps and $1/2-$Harmonic Maps into Plane Distributions.}
\label{reg}

\subsection{Variational Harmonic Maps into Plane Distributions}
\subsubsection{The $1$-D Case}
In this subsection we consider the well known case of critical points of the Dirichlet energy within the space
\[
{\frak H}^{1}(S^1):=\lf\{u\in H^{1}(S^1,{\R}^m)\quad;\quad P_N(u)\frac{d u}{d\theta}=0\quad\mbox{ in }{\mathcal D}'(S^1)\rg\}
\]
We introduce the following Lagrangian defined on the Hilbert Space $H^{1}(S^1,{\R}^m)\times L^2(S^1,{\R}^m)$ 
\be
\label{IV.1}
{\mathcal L}^1(u,\xi):=\int_{S^1}\frac{\lf<\xi,P_T(u)\xi\rg>}{2}\ d\theta-\int_{S^1}\xi\cdot\frac{du}{d\theta}\ d\theta
\ee
A point $(u,\xi)$ is a critical point to ${\mathcal L}$ if and only if for any $(w,\eta)\in H^{1}(S^1,{\R}^m)\times L^2(S^1,{\R}^m)$ we have
\be
\label{IV.2}
\begin{array}{l}
\ds\int_{S^1}\lf<\eta,P_T(u)\,\xi\rg> d\theta-\int_{S^1}\eta\cdot\frac{du}{d\theta}\ d\theta+\int_{S^1}\frac{\lf<\xi,d_wP_T(u)\xi\rg>}{2}\ d\theta-\int_{S^1}\xi\cdot\frac{dw}{d\theta}\ d\theta=0\\[5mm]

\end{array}
\ee
This is equivalent to
\be
\label{IV.2}
\lf\{
\begin{array}{l}
\ds \frac{du}{d\theta}=P_T(u)\,\xi\\[5mm]
\ds \frac{d\xi_k}{d\theta}=-\frac{1}{2}\lf<\xi,\p_{z_k}P_T(u)\xi\rg>\quad\quad\forall\,k=1\cdots m
\end{array}
\rg.
\ee
This implies first that $d\xi/d\theta\in L^1(S^1)$ which gives that $\xi\in C^0(S^1)$. Hence we deduce that $(u,\xi)\in C^1(S^1)\times C^1(S^1)$. This is the case 
of normal geodesics in sub-riemannian geometric.

\medskip

Assume that $P_T$ is \underbar{integrable}, i.e. satisfies (\ref{I.2}), then taking the $\theta$ derivative of the first equation of (\ref{IV.2}) gives
\be
\label{IV.3}
P_T(u)\, \frac{d^2u}{d\theta^2}=P_T(u)\,d_{\frac{du}{d\theta}}P_T(u)\,\xi+P_T(u)\,\frac{d\xi}{d\theta}
\ee
We have using the second equation of (\ref{IV.2})
\be
\label{IV.4}
\begin{array}{l}
\ds\lf<\ep_l,P_T(u)\,\frac{d\xi}{d\theta}\rg>=\sum_{k=1}^mP_T^{lk}\,\frac{d\xi_k}{d\theta}
=-\frac{1}{2}\sum_{i,j,k=1}^m P_T^{lk}\,\p_{z_k}P_T^{ij}\, \xi_i\,\xi_j\\[5mm]
\quad=-\frac{1}{2}\sum_{i,j,k,s=1}^m P_T^{lk}\,\p_{z_k}P_T^{is}\,P_T^{sj}\, \xi_i\,\xi_j-\frac{1}{2}\sum_{i,j,k,s=1}^m P_T^{lk}\,P_T^{is}\,\p_{z_k}P_T^{sj}\, \xi_i\,\xi_j\\[5mm]
\ds\quad=- \lf<\xi,d_{P_T\ep_l}P_T\ P_T\xi\rg>
\end{array}
\ee
Combining the previous with Lemma~\ref{lma-A.1} gives
\[
\lf<\ep_l,P_T(u)\,\frac{d\xi}{d\theta}\rg>=- \lf<\xi,d_{P_T\xi}P_T\ P_T\ep_l\rg>\\[5mm]
\]
Using the symmetry of the matrix $d_{P_T\xi}P_T$ and $P_T$ we have
\[
\lf<\ep_l,P_T(u)\,\frac{d\xi}{d\theta}\rg>=- \lf<d_{P_T\xi}P_T\,\xi,P_T\ep_l\rg>=-\lf<\ep_l,P_T\, d_{\frac{d u}{d\theta}}P_T\,\xi\rg>
\]
so in other words we have proved in the integrable case
\be
\label{IV.5}
P_T(u)\,\frac{d\xi}{d\theta}=-P_T(u)\, d_{\frac{d u}{d\theta}}P_T(u)\,\xi\quad.
\ee
Combining (\ref{IV.3}) and (\ref{IV.5}) we obtain
\[
P_T(u)\,\frac{d^2u}{d\theta^2}=0\quad.
\]
which is the well known harmonic map equation (\ref{I.4}) for $\al=1$.
\subsection{2-Dimensional Variational Harmonic Maps into Plane Distributions.}

Following the 1-dimensional case one can introduce for pairs $(u,\xi)\in W^{1,2}(D^2,{\R}^m)\times L^2(D^2,{\R}^2\otimes{\R}^m)$
\[
{\mathcal L}(u,\xi):=\int_{D^2}\frac{\sum_{i=1}^2\lf<\xi^l,P_T(u)\xi^l \rg>}{2}-\sum_{i=1}^2\xi^l\cdot\p_{x_l}u\ dx^2
\]
A pair $(u,\xi)$ is a critical point of ${\mathcal L}$ if and only if
\be
\label{IV.2-b}
\lf\{
\begin{array}{l}
\ds \frac{\p u}{\p x_l}=P_T(u)\,\xi^l\\[5mm]
\ds \mbox{div}\,\xi_k=-\frac{1}{2}\sum_{l=1}^2\lf<\xi^l,\p_{z_k}P_T(u)\,\xi^l\rg>\quad\quad\forall\,k=1\cdots m
\end{array}
\rg.
\ee
where
\[
\mbox{div}\,\xi_k=\p_{x_1}\xi^1_k+\p_{x_2}\xi^2_k\quad.
\]
Similarly as (\ref{IV.4}) we have
\[
\lf<\ep_i, P_T(u)\,\mbox{div}\,\xi\rg>=-\sum_{l=1}^2\lf<\xi^l,\p_{P_T\ep_i}P_T(u)\,P_T(u)\,\xi^l\rg>\quad\quad\forall\,k=1\cdots m\quad
\]
Hence
\[
\begin{array}{l}
\ds\mbox{div}\,(P_T(u)\xi)_i=\sum_{l=1}^2\sum_{j=1}^md_{\p_{x_l}u}P^{ij}_T(u)\,\xi^l_j-\sum_{l=1}^2\sum_{j,k=1}^m\xi^l_j\ d_{P_T\ep_i}P^{jk}_T(u)\, (P_T(u)\xi^l)_k\\[5mm]
\ds\quad=\sum_{l=1}^2\sum_{j,k=1}^m\xi^l_j\,\p_{z_k}P^{ij}_T(u)\, (P_T(u)\xi^l)_k-\sum_{l=1}^2\sum_{j,k,s=1}^m\xi^l_j\ P^{is}_T\,\p_{z_s}P^{jk}_T(u)\, (P_T(u)\xi^l)_k\\[5mm]
\ds\quad=\sum_{l=1}^2\sum_{j,k,s=1}^m\xi^l_j\,P^{ks}_T(u)\,\p_{z_s}P^{ij}_T(u)\, (P_T(u)\xi^l)_k-\sum_{l=1}^2\sum_{j,k,s=1}^m\xi^l_j\ P^{is}_T(u)\,\p_{z_s}P^{kj}_T(u)\, (P_T(u)\xi^l)_k
\end{array}
\]
Denote
\[
\Om^{ik}_l:=\sum_{j,s=1}^m\xi^l_j\,P^{is}_T\,\p_{z_s}P^{kj}_T(u)-\xi^l_j\ P^{ks}_T\,\p_{z_s}P^{ij}_T(u)\,
\]
We have by definition
\[
\forall\, l=1,2\quad\forall\,i,k\in\{1\cdots m\}\quad\quad\Om^{ik}_l=-\,\Om^{ki}_l
\]
Moreover $\Om\in L^2$ and $u$ satisfies the following system
\[
-\,-\Delta u=\Omega\cdot\nabla u\quad\mbox{ in }D^2
\]
Hence, using \cite{Riv}, we have that $\nabla u\in \cap_{p<2} W^{1,p}_{loc}(D^2)$.

\subsection{Variational $1/2-$Harmonic Maps into Plane Distributions}
On $ H^{s}(S^1)$ ($s\in{\R}$ arbitrary) we define the following operator for any $\al\in [0,1]$
\[
(-\Delta)^{-\al/2}_0\ \quad f\in H^s(S^1)\longrightarrow v=(-\Delta)^{-\al/2}_0f\in H^{s+\al}_0(S^1)
\]
where $v$ satisfies
\[
(-\Delta)^{\al/2}v=f-\frac{1}{2\pi}\int_{S^1}f
\]
and is given explicitly by
\[
v:=\sum_{n\in {\Z}^\ast} f_n |n|^{-\al}\ e^{in\theta}\quad\mbox{where}\quad f=\sum_{n\in {\Z}}f_n\, e^{i\,n\,\theta}\quad.
\]
Observe that $v$ satisfies
\[
\int_{S^1}v(\theta)\ d\theta=0\quad.
\]
Hence this gives in particular that
\be
\label{IV.5-a}
(-\Delta)^{-\al/2}_0\circ(-\Delta)^{-\beta/2}_0 f=(-\Delta)^{-(\al+\beta)/2}_0\, f
\ee
We have also for any $f\in H^{-\al}(S^1,{\R}^m)$ and $g\in L^2(S^1,{\R}^m)$
\be
\label{IV.5-b}
\int_{S^1}(-\Delta)_0^{-\al/2}f(\theta)\  g(\theta)\ d\theta=(f,(-\Delta)_0^{-\al/2}\, g)_{H^{-\al},H^\al}
\ee
We introduce the following Lagrangian defined on the sub-manifold of the Hilbert Space $H^{1/2}(S^1,{\R}^m)\times H^{-1/2}(S^1,{\R}^m)$
given by
\[
{\mathfrak E}:=\lf\{
\begin{array}{c}
\ds(u,\xi)\in H^{1/2}(S^1,{\R}^m)\times H^{-1/2}(S^1,{\R}^m)\quad\mbox{s. t. }\\[5mm]
\ds\lf(P_N(u),\frac{du}{d\theta}\rg)_{H^{1/2},H^{-1/2}}=0\\[5mm]
\ds(-\Delta)_0^{-1/4}(P_T(u)\xi)\in L^2(S^1)\quad\mbox{ and }\quad(-\Delta)^{-1/4}_0\lf(P_T(u)\frac{du}{d\theta}\rg)\in L^2(S^1)\quad 
\end{array}\rg\}
\]
Let
\[
\begin{array}{l}
\ds{\mathcal L}^{1/2}(u,\xi):=\int_{S^1}\frac{|(-\Delta)^{-1/4}_0(P_T(u)\xi)|^2}{2}\ d\theta-\\[5mm]
\ds \quad\quad\quad\int_{S^1}\lf<(-\Delta)^{-1/4}_0(P_T(u)\xi),(-\Delta)^{-1/4}_0\lf(P_T(u)\frac{du}{d\theta}\rg)\rg>\ d\theta\\[5mm]
\ds \quad\quad\quad-\int_{S^1}\lf<(-\Delta)^{-1/4}_0(P_N(u)\xi),(-\Delta)^{-1/4}_0\lf(P_N(u)\frac{du}{d\theta}\rg)\rg>\ d\theta
\end{array}
\]
Observe that if $u\in {\frak H}^{1/2}$ we have $(u,du/d\theta)\in {\mathfrak E}$ and
\[
{\mathcal L}^{1/2}\lf(u,\frac{du}{d\theta}\rg)=-\frac{1}{2}\int_{S^1}|(-\Delta)^{1/4}u|^2\ d\theta\quad.
\]
Assume now $(u,\xi)$ is a critical point of ${\mathcal L}^{1/2}$ in ${\mathfrak E}$. Hence for any choice of $(w,\eta)\in C^\infty(S^1,{\R}^m)\times  C^\infty(S^1,{\R}^m)$ where $w$ satisfies the constraint
\be
\label{IV.5c}
\int_{S^1}P_N(u)\,\frac{dw}{d\theta}\, d\theta+\lf<d_wP_N(u)\,,\,\frac{du}{d\theta}\rg>_{H^{1/2},{H^{-1/2}}}=0
\ee
we have respectively
\be
\label{IV.6}
\begin{array}{l}
\ds\int_{S^1}\lf<(-\Delta)^{-1/4}_0(P_T(u)\xi),(-\Delta)^{-1/4}_0(P_T(u)\eta)   \rg>\ d\theta\\[5mm]
\ds\quad\quad\quad-\int_{S^1}\lf<(-\Delta)^{-1/4}_0(P_T(u)\eta),(-\Delta)^{-1/4}_0\lf(P_T(u)\frac{du}{d\theta}\rg)\rg>\ d\theta\\[5mm]
\ds \quad\quad\quad-\int_{S^1}\lf<(-\Delta)^{-1/4}_0(P_N(u)\eta),(-\Delta)^{-1/4}_0\lf(P_N(u)\frac{du}{d\theta}\rg)\rg>\ d\theta=0
\end{array}
\ee
and
\be
\label{IV.6}
\begin{array}{l}
\ds\int_{S^1}\lf<(-\Delta)^{-1/4}_0(P_T(u)\xi),(-\Delta)^{-1/4}_0(d_wP_T(u)\xi)   \rg>\ d\theta\\[5mm]
\ds \quad\quad\quad-\int_{S^1}\lf<(-\Delta)^{-1/4}_0(d_wP_T(u)\xi),(-\Delta)^{-1/4}_0\lf(P_T(u)\frac{du}{d\theta}\rg)\rg>\ d\theta\\[5mm]
\ds\quad\quad\quad-\int_{S^1}\lf<(-\Delta)^{-1/4}_0(P_T(u)\xi),(-\Delta)^{-1/4}_0\lf(d_wP_T(u)\frac{du}{d\theta}\rg)\rg>\ d\theta\\[5mm]
\ds \quad\quad\quad-\int_{S^1}\lf<(-\Delta)^{-1/4}_0(d_wP_N(u)\xi),(-\Delta)^{-1/4}_0\lf(P_N(u)\frac{du}{d\theta}\rg)\rg>\ d\theta\\[5mm]
\ds \quad\quad\quad-\int_{S^1}\lf<(-\Delta)^{-1/4}_0(P_N(u)\xi),(-\Delta)^{-1/4}_0\lf(d_wP_N(u)\frac{du}{d\theta}\rg)\rg>\ d\theta\\[5mm]
\ds \quad\quad\quad-\int_{S^1}\lf<(-\Delta)^{-1/4}_0(P_T(u)\xi),(-\Delta)^{-1/4}_0\lf(P_T(u)\frac{dw}{d\theta}\rg)\rg>\ d\theta\\[5mm]
\ds \quad\quad\quad-\int_{S^1}\lf<(-\Delta)^{-1/4}_0(P_N(u)\xi),(-\Delta)^{-1/4}_0\lf(P_N(u)\frac{dw}{d\theta}\rg)\rg>\ d\theta=0
\end{array}
\ee
The first equation (\ref{IV.6}) implies using (\ref{IV.5-a}), (\ref{IV.5-b}) and the symmetry of the matrices $P_T$ and $P_N$
\be
\label{IV.7}
P_T(u)\,(-\Delta_0)^{-1/2}\lf(P_T(u)\,\xi-P_T(u)\frac{du}{d\theta}\rg)-P_N(u)\,(-\Delta_0)^{-1/2}\lf(P_N(u)\frac{du}{d\theta}\rg)=0
\ee
This implies
\be
\label{IV.8}
\lf\{
\begin{array}{l}
\ds P_T(u)\, (-\Delta_0)^{-1/2}\lf(P_T(u)\,\xi-P_T(u)\frac{du}{d\theta}\rg)=0\\[5mm]
\ds P_N(u)\, (-\Delta_0)^{-1/2}\lf(P_N(u)\frac{du}{d\theta}\rg)=0
\end{array}
\rg.
\ee
Multiplying the second equation by $du/d\theta$ and integrating by parts gives
\[
\int_{S^1}\lf|(-\Delta)^{-1/4}_0\lf(P_N(u)\frac{du}{d\theta}\rg)\rg|\ d\theta=0
\]
this gives
\be(-\Delta)^{-1/2}
\label{IV.9}
P_N(u)\frac{du}{d\theta}\equiv Cte
\ee
Since the membership of $u$ to ${\mathfrak E}$ imposes
\[
\lf(P_N(u),\frac{du}{d\theta}\rg)_{H^{1/2},H^{-1/2}}=0
\]
Hence we have
\be
\label{IV.10}
P_N(u)\frac{du}{d\theta}\equiv 0
\ee
or in other words $u\in {\mathcal H}$. Multiplying now the first equation by $\xi-du/d\theta$ and integrating by parts gives
\be
\label{IV.11}
P_T(u)\,\xi-\frac{du}{d\theta}\equiv Cte
\ee
With these informations at hand (\ref{IV.6}) becomes
\be
\label{IV.12}
\begin{array}{l}
\ds-\int_{S^1}\lf<\frac{dw}{d\theta},  P_T(u)\,(-\Delta)^{-1/2}_0\lf(P_T(u)\,\xi\rg)+ P_N(u)\,(-\Delta)^{-1/2}_0\lf(P_N(u)\,\xi\rg)\rg>\\[5mm]
\ds+\int_{S^1}\lf<d_wP_T(u)\,\frac{du}{d\theta},(-\Delta_0)^{-1/2}(P_N(u)\,\xi)\rg>\ d\theta\\[5mm]
\ds-\int_{S^1}\lf<d_wP_T(u)\,\frac{du}{d\theta},(-\Delta_0)^{-1/2}(P_T(u)\,\xi)\rg>\ d\theta=0
\end{array}
\ee
Combining (\ref{IV.5c}) and (\ref{IV.12}) and assuming $u$ is a non degenerate point of the constraint
\[
\lf(P_N(u),\frac{du}{d\theta}\rg)_{H^{1/2},H^{-1/2}}=0
\]
we obtain the existence of $\la=(\la_1\cdots\la_m)\in{\R}^m$ such that for any $k=1\cdots m$
\be
\label{IV.13}
\begin{array}{l}
\ds\frac{d}{d\theta}\lf(P_T(u)\,(-\Delta)^{-1/2}_0\lf(P_T(u)\,\xi\rg)+ P_N(u)\,(-\Delta)^{-1/2}_0\lf(P_N(u)\,\xi\rg)\rg)^k=\\[5mm]
\ds\quad-\lf<\p_{z_k}P_T(u)\,\frac{du}{d\theta},(-\Delta_0)^{-1/2}(P_N(u)\,\xi)\rg>+\lf<\p_{z_k}P_T(u)\,\frac{du}{d\theta},(-\Delta_0)^{-1/2}(P_T(u)\,\xi)\rg>\\[5mm]
\ds\quad+\lf<\la,\p_{z_k}P_T(u)\,\frac{du}{d\theta}-\p_{\frac{du}{d\theta}}P_T(u)\,\ep_k\rg>
\end{array}
\ee

\medskip

Assume that $P_T$ is \underbar{integrable}, i.e. satisfies (\ref{I.2}). Taking the multiplication of (\ref{IV.13}) with $P_T(u)$ gives
\be
\label{IV.14}
\begin{array}{l}
\ds  \lf<\ep_i, P_T(u)\, {\mathcal R}\, P_T(u)\,\xi  \rg>\\[5mm]
\ds  +\lf<\ep_i, P_T(u)\,\frac{d P_T(u)}{d\theta}\, (-\Delta)^{-1/2}_0\lf(P_T(u)\,\xi\rg) +P_T(u)\,\frac{d P_N(u)}{d\theta}\, (-\Delta)^{-1/2}_0\lf(P_N(u)\,\xi\rg)  \rg>\\[5mm]
\ds=-\lf<\p_{P_T(u)\,\ep_i}P_T(u)\,\frac{du}{d\theta},(-\Delta_0)^{-1/2}(P_N(u)\,\xi)\rg>\\[5mm]
\ds\quad+\lf<\p_{P_T(u)\,\ep_i}P_T(u)\,\frac{du}{d\theta},(-\Delta_0)^{-1/2}(P_T(u)\,\xi)\rg>\\[5mm]
\ds\quad+\lf<\la,\p_{P_T(u)\,\ep_i}P_T(u)\,\frac{du}{d\theta}-\p_{\frac{du}{d\theta}}P_T(u)\,P_T(u)\,\ep_i\rg>

\end{array}
\ee
where ${\mathcal R}$ is the Riesz operator given by
\[
{\mathcal R}\ :\ f=\sum_{n\in {\Z}}f_n\, e^{i\,n\,\theta}\ \longrightarrow\ {\mathcal R}f:=i\,\sum_{n\in {\Z}^\ast} \mbox{sgn}(n) \,f_n\, e^{i\,n\,\theta}
\]
Since $P_N(u)\frac{du}{d\theta}=0$, on can use lemma~\ref{lma-A.1} in order to infer
\be
\label{IV.15}
\lf<\la,\p_{P_T(u)\,\ep_i}P_T(u)\,\frac{du}{d\theta}-\p_{\frac{du}{d\theta}}P_T(u)\,P_T(u)\,\ep_i\rg>=0
\ee
moreover, using again lemma~\ref{lma-A.1}, the symmetry of the matrices $d P_N(u)/d\theta$ and $P_T(u)$, we obtain
\be
\label{IV.16}
\begin{array}{l}
\lf<\p_{P_T(u)\,\ep_i}P_T(u)\,\frac{du}{d\theta},(-\Delta_0)^{-1/2}(P_N(u)\,\xi)\rg>=\lf<\p_{\frac{du}{d\theta}}P_T(u)\,P_T(u)\,\ep_i,(-\Delta_0)^{-1/2}(P_N(u)\,\xi)\rg>\\[5mm]
\ds\quad=-\lf<\p_{\frac{du}{d\theta}}P_N(u)\,P_T(u)\,\ep_i,(-\Delta_0)^{-1/2}(P_N(u)\,\xi)\rg>\\[5mm]
\ds\quad=-\lf<\frac{d P_N(u)}{d\theta}\,P_T(u)\,\ep_i,(-\Delta_0)^{-1/2}(P_N(u)\,\xi)\rg>\\[5mm]
\ds\quad=-\lf<P_T(u)\,\ep_i,\frac{d P_N(u)}{d\theta}\,(-\Delta_0)^{-1/2}(P_N(u)\,\xi)\rg>\\[5mm]
\ds\quad=-\lf<\ep_i,P_T(u)\frac{d P_N(u)}{d\theta}\,(-\Delta_0)^{-1/2}(P_N(u)\,\xi)\rg>
\end{array}
\ee
and similarly we have
\be
\label{IV.17}
\begin{array}{l}
\lf<\p_{P_T(u)\,\ep_i}P_T(u)\,\frac{du}{d\theta},(-\Delta_0)^{-1/2}(P_T(u)\,\xi)\rg>=\lf<\p_{\frac{du}{d\theta}}P_T(u)\,P_T(u)\,\ep_i,(-\Delta_0)^{-1/2}(P_T(u)\,\xi)\rg>\\[5mm]
\ds\quad=\lf<\ep_i,P_T(u)\frac{d P_T(u)}{d\theta}\,(-\Delta_0)^{-1/2}(P_T(u)\,\xi)\rg>
\end{array}
\ee
Combining (\ref{IV.14})...(\ref{IV.17}) we obtain
\be
\label{IV.18}
\begin{array}{l}
 0=P_T(u)\, {\mathcal R}\, P_T(u)\,\xi=P_T(u)\, {\mathcal R}\,\frac{du}{d\theta}=P_T(u)\,(-\Delta)^{1/2}u
 \end{array}
\ee
which is exactly the $1/2-$harmonic map equation.

\medskip

In fact the correspondence between critical points of ${\mathcal L}^{1/2}$ and critical points of the $1/2-$energy within ${\frak H}^{1/2}$ goes beyond the very
special case of integrable plane distributions. Precisely we have the following theorem.
\begin{Th}
\label{th-IV.1}
Let $(u,\xi)$ be a smooth critical point of ${\mathcal L}^{1/2}$ in ${\mathfrak E}$ then $u$ is a critical point of 
\[
E^{1/2}(u)=\int_{S^1}|(-\Delta)^{1/4}u|^2 \ d\theta
\]
within the space ${\frak H}^{1/2}$ of horizontal $H^{1/2}-$maps.\hfill $\Box$
\end{Th} 
\noindent{\bf Proof of theorem~\ref{th-IV.1}.}
Since $u$ is assumed to be smooth we can make use locally of an orthonormal frame ${e}_1\cdots {e}_n$ generating the plane distribution
given by the Images of $P_T$. With this frame at hand we can introduce the {\it control} $\al_1(\theta),\cdots,\al_n(\theta)$ such that
\be
\label{IV.18-a}
\frac{d u}{d\theta}=\sum_{i=1}^n\al_i(\theta)\,{e}_i(u(\theta))
\ee
Classical considerations from control theory in sub-riemannian framework (see for instance \cite{Mon}) asserts that an infinitesimal variation of an horizontal map satisfying (\ref{IV.18-a})
is given by $w$ satisfying
\be
\label{IV.18-b}
\frac{d w}{ d\theta}=\sum_{i=1}^n v_i(\theta)\ e_i(u(\theta))+\sum_{i=1}^n\al_i(\theta)\ d_w e_i(u(\theta))\quad.
\ee
where the $v_i(\theta)$ are arbitrary so that the constraint (\ref{IV.5c}) is satisfied. Since $P_T\,P_N=0$, we have
\[
d_wP_T(u)\, P_N+P_T(u)\,d_wP_N(u)=0
\]
Hence this implies, using that $d_w P_T=- \,d_w P_N$,
\be
\label{IV.18-c}
P_T(u)\, d_wP_T(u)\, P_T(u)=-P_T(u)\, d_wP_N(u)\, P_T(u)=d_wP_T(u)\, P_N(u)\,P_T(u)=0
\ee
Hence
\be
\label{IV.18-d}
\begin{array}{l}
\ds P_T(u)\,d_wP_T(u)\,\frac{du}{d\theta}=P_T(u)\, d_wP_T(u)\, P_T(u)\,\frac{du}{d\theta}=0\\[5mm]
\ds\quad\quad\quad\Longrightarrow\quad d_wP_T(u)\,\frac{du}{d\theta}=P_N(u)\,d_wP_T(u)\,\frac{du}{d\theta}\quad.
\end{array}
\ee
Since 
\[
P_T:=\sum_{i=1}^ne_i\otimes e_i
\]
We have that
\be
\label{IV.18-e}
d_wP_T(u)\,\frac{du}{d\theta}=\sum_{i,j=1}^n\al_j\ d_we_i\cdot e_j\ e_i+\sum_{i=1}^n\al_i\ d_we_i\quad.
\ee
Combining (\ref{IV.18-b}), (\ref{IV.18-d}) and (\ref{IV.18-e}) we obtain
\be
\label{IV.18-f}
d_wP_T(u)\,\frac{du}{d\theta}=\sum_{i=1}^n\al_i\ P_N(u)\,d_we_i=P_N\,\frac{d w}{d\theta}\quad.
\ee
Inserting this identity in (\ref{IV.12})
\be
\label{IV.18-g}
\begin{array}{l}
\ds-\int_{S^1}\lf<\frac{dw}{d\theta},  P_T(u)\,(-\Delta)^{-1/2}_0\lf(P_T(u)\,\xi\rg)+ P_N(u)\,(-\Delta)^{-1/2}_0\lf(P_N(u)\,\xi\rg)\rg>\ d\theta\\[5mm]
\ds+\int_{S^1}\lf<P_N\,\frac{d w}{d\theta},(-\Delta_0)^{-1/2}(P_N(u)\,\xi)\rg>\ d\theta\\[5mm]
\ds-\int_{S^1}\lf<P_N\,\frac{d w}{d\theta},(-\Delta_0)^{-1/2}(P_T(u)\,\xi)\rg>\ d\theta=0
\end{array}
\ee
which is equivalent to
\be
\label{IV.18-h}
\int_{S^1}\lf<\frac{dw}{d\theta},  (-\Delta)^{-1/2}_0\lf(\frac{du}{d\theta}\rg)\rg>\ d\theta=0
\ee
Since this holds for any perturbation $w$ of $u$ in ${\frak H}^{1/2}$, we have proved the theorem.\hfill $\Box$

\subsection{Reformulation of  the Euler-Lagrange Equation}
Observe that (\ref{IV.13}) becomes
\be
\label{IV.22}
\begin{array}{l}
\ds \lf(\lf(P_T(u)\,{\mathcal R}\,P_T(u)+P_N(u)\,{\mathcal R}\,P_N(u)\rg)\xi\rg)^k\\[5mm]
\ds+\lf(\frac{d P_T(u)}{d\theta}\lf((-\Delta)^{-1/2}_0\lf(P_T(u)\,\xi\rg)\rg)+ \frac{d P_N(u)}{d\theta}\,(-\Delta)^{-1/2}_0\lf(P_N(u)\,\xi\rg)\rg)^k=\\[5mm]
\ds\quad+\lf<\p_{z_k}P_N(u)\,\frac{du}{d\theta},(-\Delta_0)^{-1/2}(P_N(u)\,\xi)\rg>+\lf<\p_{z_k}P_T(u)\,\frac{du}{d\theta},(-\Delta_0)^{-1/2}(P_T(u)\,\xi)\rg>\\[5mm]
\ds\quad+\lf<\la,\p_{z_k}P_T(u)\,\frac{du}{d\theta}-\p_{\frac{du}{d\theta}}P_T(u)\,\ep_k\rg>
\end{array}
\ee
This gives
\be
\label{IV.23}
\begin{array}{l}
\ds \lf(\lf(P_T(u)\,{\mathcal R}\,P_T(u)+P_N(u)\,{\mathcal R}\,P_N(u)\rg)\xi\rg)^k=\\[5mm]
\ds\sum_{j=1}^m \lf(\sum_{i=1}^m(\p_{z_k}P_N^{ij}-\p_{z_j}P_N^{ik})\ (-\Delta_0)^{-1/2}(P_N(u)\,\xi)^i+(\p_{z_k}P_T^{ij}-\p_{z_j}P_T^{ik})\ (-\Delta_0)^{-1/2}(P_T(u)\,\xi)^i\rg)\, \frac{du^j}{d\theta}\\[5mm]
\ds+\sum_{j=1}^m \lf(\sum_{i=1}^m(\p_{z_k}P_T^{ij}-\p_{z_j}P_T^{ik})\, \la^i\rg)\, \frac{du^j}{d\theta}
\end{array}
\ee
Denote
\[
\begin{array}{l}
\ds\om^{kj}:=\lf(\sum_{i=1}^m(\p_{z_k}P_N^{ij}-\p_{z_j}P_N^{ik})\ (-\Delta_0)^{-1/2}(P_N(u)\,\xi)^i+(\p_{z_k}P_T^{ij}-\p_{z_j}P_T^{ik})\ ((-\Delta_0)^{-1/2}(P_T(u)\,\xi)^i+\la^i)\rg)\\[5mm]
\ds\quad\quad=\lf(\sum_{i=1}^m(\p_{z_k}P_N^{ij}-\p_{z_j}P_N^{ik})\ (-\Delta_0)^{-1/2}\xi^i+(\p_{z_k}P_T^{ij}-\p_{z_j}P_T^{ik})\ \la^i\rg)
\end{array}
\]
Observe that $\om$ is \underbar{antisymmetric} and the equation becomes
\be
\label{IV.24}
\lf(P_T\,{\mathcal R}\,P_T+P_N\,{\mathcal R}\,P_N\rg)\xi=\om\,P_T\, \xi
\ee 
Let
\[
v=\lf(\begin{array}{c} P_T\,\xi\\[3mm]
{\mathcal  R}\,P_N\,\xi
\end{array}
\rg)
\]
Observe that
\be
\label{IV.25}
(P_T{\mathcal R}+P_N)\, v=\lf(\begin{array}{c}P_T\, {\mathcal R}\,P_T\,\xi\\[3mm]
P_N\,{\mathcal  R}\,P_N\,\xi
\end{array}
\rg)
\ee
If one multiplies (\ref{IV.24}) by $P_T(u)$ one gets that $w:=P_T(u)\xi$ satisfies
\be
\label{IV.26}
\lf\{
\begin{array}{l}
\ds  P_T\,{\mathcal R}\, w=\Om\, w\\[5mm]
\ds P_N \,w=0
\end{array}
\rg.
\ee
where $\Om:=P_T\,\om\,P_T$ is \underbar{antisymmetric} which is a ``deformation'' of the $1/2-$harmonic equation
\be
\label{IV.26}
\lf\{
\begin{array}{l}
\ds  P_T\,{\mathcal R}\, w=0\\[5mm]
\ds P_N \,w=0
\end{array}
\rg.
\ee
where $w:=du/d\theta$.

\renewcommand{\theequation}{A.\arabic{equation}}
\renewcommand{\theTh}{A.\arabic{Th}}
\renewcommand{\theProp}{A.\arabic{Prop}}
\renewcommand{\theLma}{A.\arabic{Lma}}
\renewcommand{\theCo}{A.\arabic{Co}}
\renewcommand{\theRm}{A.\arabic{Rm}}
\renewcommand{\theequation}{A.\arabic{equation}}
\setcounter{equation}{0} 

\appendix

\section{Appendix}

\subsection{Integrable Distributions}
The goal of the present section is to establish the following elementary lemma which is well known.
\begin{Lma}
\label{lma-A.1}
Let $P_T$ be a $C^1$ plane distribution satisfying (\ref{I.1}) and assume $P_T$ is integrable, i.e. satisfies (\ref{I.2}),
then
\be
\label{A.1}
\forall\ X,Y\in C^1({\R}^m,{\R}^m)\quad\quad \mbox{ we have }\quad d_{P_TX}P_T\ P_TY=d_{P_TY}P_T\ P_TX
\ee
or in other words
\be
\label{A.2}
\forall\ i\,j\,k\in\{1\cdots m\}\quad\quad\sum_{s,t=1}^m\p_{z_t}P_T^{is}\ P^{sk}_T\ P^{tj}_T=\sum_{s,t=1}^m\p_{z_t}P_T^{is}\ P^{sj}_T\ P^{tk}_T
\ee
\hfill $\Box$
\end{Lma}
\noindent{\bf Proof of lemma~\ref{lma-A.1}.} Let $(\ep_i)_{i=1\cdots m}$ be the canonical basis of ${\R}^m$. We have
\[
[P_T\,\ep_j,P_T\,\ep_k]=\sum_{s,t=1}^m\lf(P^{tj}_T\, \p_{z_t}P^{sk}_T-P^{tk}_T\, \p_{z_t}P^{sj}_T\rg)\ \ep_s
\]
Equation (\ref{I.2}) becomes
\[
\forall\,i,j,k\quad\quad\sum_{s,t=1}^m(\delta^{is}-P_T^{is})\ \lf(P^{tj}_T\, \p_{z_t}P^{sk}_T-P^{tk}_T\, \p_{z_t}P^{sj}_T\rg)=0
\]
which gives
\be
\label{A.3}
\begin{array}{l}
\ds 0=\sum_{t=1}^mP^{tj}_T\, \p_{z_t}P^{ik}_T-P^{tk}_T\, \p_{z_t}P^{ij}_T-\sum_{s,t=1}^m P^{tj}_T\, P_T^{is}\,\p_{z_t}P^{sk}_T-P^{tk}_T\,P_T^{is}\, \p_{z_t}P^{sj}_T
\end{array}
\ee
Using the fact that $P_T\circ P_T=P_T$ we have
\be
\label{A.4}
-\sum_{s=1}^mP_T^{is}\,\p_{z_t}P^{sk}_T=-\p_{z_t}P^{ik}_T+\sum_{s=1}^m\p_{z_t}P_T^{is}\,P^{sk}_T
\ee
Combining (\ref{A.3}) and (\ref{A.4}) gives then (\ref{A.2}) and lemma~\ref{lma-A.1} is proved.\hfill $\Box$

\subsection{Rewriting the Commutators}\label{rewrite}
In this section we recall the  explicit form of the matrices $\Omega_0,\Omega_1$ and of the operator ${\cal{Z}}$ introduced in \rec{modelsystem} in the case of $1/2$-harmonic maps.
\begin{Proposition}\label{EulEq}
Let $u\in\dot{H}^{1/2}(\R,{\cal{N}} )$ be a weak $1/2$-harmonic map. Then the following equation holds
\begin{eqnarray}\label{Eulerbis}
\Delta^{1/4}v=(-\Delta)^{1/4}\left(\begin{array}{l}
P_T(-\Delta)^{1/4} u\\ 
{\cal{R}}P_N(-\Delta)^{1/4} u\end{array}\right)&=& \tilde \Omega+   \Omega_1\left(\begin{array}{l}
P_T(-\Delta)^{1/4} u\\ 
{\cal{R}}P_N(-\Delta)^{1/4} u\end{array}\right) \\ 
&+&\Omega
\left(\begin{array}{l}
P_T(-\Delta)^{1/4} u\\
{\cal{R}}P_N(-\Delta)^{1/4} u\end{array}\right)\,,\nonumber
\end{eqnarray}
 where $\Om=\Omega(P_T)\in L^2({\R},so(2m))$, $ \Omega_1=\Omega_1(P_T) \in L^{2,1}$ with
 $$
 \|\Om\|_{L^2},\| \Omega_1\|_{L^{2,1}}\le C(\|P_T\|_{{H}^{1/2}}+\|P_T\|^2_{{H}^{1/2}}),$$
 and
 \begin{equation}
{ \tilde\Omega=(C-2D)\left( \begin{array}{c} P_T(-\Delta)^{1/4}u\\{\cal{R}}P_N(-\Delta)^{1/4}u\end{array}\right)}
\end{equation}
where  
the matrices $C$ and $D$  are $2\times 2m$ matrices  whose components are made by pseudo-differential operators: for $j\in\{1,\ldots, m\}$
{ \begin{equation}
 \begin{array}{c}
c_{1j}=(-\Delta)^{1/4}\{P_T\}-P_T(-\Delta)^{1/4}+(-\Delta)^{1/4}[P_T]\\[5mm]
c_{1,j+m}=(-\Delta)^{1/4}\{P_T\circ {\cal{R}}\}-P_T(-\Delta)^{1/4}\circ {\cal{R}}+(-\Delta)^{1/4}[P_T]\circ {\cal{R}}\\[5mm]
c_{2,j}={\cal{R}}(-\Delta)^{1/4}\{P_N\}-P_N(-\Delta)^{1/4}\circ {\cal{R}}-(-\Delta)^{1/4}[P_N]\circ {\cal{R}}\\[5mm]
c_{2,j+m}=(-\Delta)^{1/4}\{P_N\}+P_N(-\Delta)^{1/4} -(-\Delta)^{1/4}[P_N]\,. 
\end{array}
\end{equation}}
and
{ 
\begin{equation}
 \begin{array}{c}
 d_{1,j}=0\\
 d_{1,j+m}={\mathcal{R}}[\omega_1]+\omega_1{\mathcal{R}}\\[5mm]
 d_{2,j}=(-\Delta)^{1/4} P_N+{\cal{R}}((-\Delta)^{1/4} P_N) {\cal{R}}\\[5mm]
 d_{2,j+m}=[{\mathcal{R}}[\omega_2]+\omega_2{\mathcal{R}}]-((-\Delta)^{1/4} P_N){\cal{R}}-{\cal{R}}((-\Delta)^{1/4} P_N).
 \end{array}
 \end{equation}
 Moreover for $w\in L^2$,
 \begin{equation}\label{bij}
\|c_{ij}(w)\|_{  {\cal{H}}^{1}(\R)},\|d_{ij}(w)\|_{  {\cal{H}}^{1}(\R)}\le C [\|P_T\|_{\dot H^{1/2}(\R))}+\|P_T\|^2_{\dot H^{1/2}(\R))}]\|w\|_{L^{2}}\,.
\end{equation}}
\end{Proposition}
 
 In order to prove Proposition \ref{EulEq} we recall the following Proposition (Proposition 1.1 in \cite{DLR2}).

\begin{Proposition}\label{EulEqold}
Let $u\in\dot{H}^{1/2}(\R,{\cal{N}} )$ be a weak $1/2$-harmonic map. Then the following equation holds
\begin{eqnarray}\label{Eulerbis}
\Delta^{1/4}v=(-\Delta)^{1/4}\left(\begin{array}{l}
P_T(-\Delta)^{1/4} u\\ 
{\cal{R}}P_N(-\Delta)^{1/4} u\end{array}\right)&=& \tilde \Omega+   \Omega_1\left(\begin{array}{l}
P_T(-\Delta)^{1/4} u\\ 
{\cal{R}}P_N(-\Delta)^{1/4} u\end{array}\right) \\ 
&+&\Omega
\left(\begin{array}{l}
P_T(-\Delta)^{1/4} u\\
{\cal{R}}P_N(-\Delta)^{1/4} u\end{array}\right)\,,\nonumber
\end{eqnarray}
 where $\Om=\Omega(P_T)\in L^2({\R},so(2m))$, $ \Omega_1=\Omega_1(P_T) \in L^{2,1}$ with
 $$
 \|\Om\|_{L^2},\| \Omega_1\|_{L^{2,1}}\le C(\|P_T\|_{{H}^{1/2}}+\|P_T\|^2_{{H}^{1/2}}),
 $$
 $$
\tilde \Omega=\left(\begin{array}{c}
- 2F(  \omega_1(P_T) ,(P_N \Delta ^{1/4} u))+ T(P_T,u)\\[5mm]
-2F({\cal{R}}((-\Delta)^{1/4} P_N),{\cal{R}}((-\Delta)^{1/4} u))-2F(\omega_{2}(P_T), P_N((-\Delta)^{1/4} u)+{\cal{R}}(S(P_N,u))
\end{array}\right)\,$$
$ \omega_1(P_T), \omega_{2}(P_T)\in L^2$ and
$$\| \omega_{1}P_T)\|_{L^2},~\|\omega_{2}(P_T)\|_{L^2}\le C(\|P_T\|_{{H}^{1/2}}+\|P_T\|^2_{{H}^{1/2}}).$$
\end{Proposition}
 
 \vskip2cm

{\bf Proof of Proposition \ref{EulEq}.}
{We next rewrtite the matrix $\tilde\Omega$ as the product of a matrix of pseudodifferential operators times $v=(P_T(-\Delta)^{1/4}u,{\cal{R}}P_N(-\Delta)^{1/4}u)^t.$\par
{\bf Step 1.} 
We   rewrite $T(P_T,(-\Delta)^{1/4}u)$ and ${\cal{R}}(S(P_N,(-\Delta)^{1/4}u))$ in terms of $P_T(-\Delta)^{1/4}u$ and  ${\cal{R}} P_N(-\Delta)^{1/4}u$.
We observe that 
$$(-\Delta)^{1/4}u=P_N(-\Delta)^{1/4}u- {\cal{R}}[{\cal{R}}P_N(-\Delta)^{1/4}u)].$$
By linearity we get
\begin{equation}
T(P_T,  (-\Delta)^{1/4}u)=T_1(P_T, P_T(-\Delta)^{1/4}u)+T_2(P_T,  {\cal{R}}P_N(-\Delta)^{1/4}u).
\end{equation}
where for $Q\in H^{1/2}  $,  $v_1,v_2\in L^2$ it holds
\begin{eqnarray}\label{t1}
T_1(Q, v_1)&:=&(-\Delta)^{1/4}[Qv_1]\\
&&-Q(-\Delta)^{1/4}[v_1]+(-\Delta)^{1/4}[Q]  v_1.\nonumber
\end{eqnarray}

 \begin{eqnarray}\label{t2}
T_2(Q, v_2)&:=&-(-\Delta)^{1/4}[Q {\cal{R}}[v_2]]\\
&&+Q(-\Delta)^{1/4}[ {\cal{R}}v_2]-(-\Delta)^{1/4} [Q] {\cal{R}}v_2.\nonumber
\end{eqnarray}

Moreover
\begin{equation}
 {\cal{R}}[S(P_N,  (-\Delta)^{1/4}u)]=S_1(P_N, P_T(-\Delta)^{1/4}u)+S_2(P_N,  {\cal{R}}P_N(-\Delta)^{1/4}u).
\end{equation}

\begin{eqnarray}\label{s1}
S_1(Q, v_1)&:=& {\cal{R}}(-\Delta)^{1/4}[Qv_1]\\
&&Q(-\Delta)^{1/4}{\cal{R}} [v_1]-(-\Delta)^{1/4}[Q] {\cal{R}} [v_1].\nonumber
\end{eqnarray}

 \begin{eqnarray}\label{s2}
S_2(Q, v_2)&:=&  (-\Delta)^{1/4}[Q  v_2]\\
&&+Q(-\Delta)^{1/4} [  v_2]-(-\Delta)^{1/4}[Q] v_2.\nonumber
\end{eqnarray}
{ {We introduce some notations: for $Q\in H^{1/2}\cap L^\infty$  we denote by 
$$(-\Delta)^{1/4}\{Q\},~~{\cal{R}}\circ (-\Delta)^{1/4}\{Q\},~~(-\Delta)^{1/4}\{Q\circ {\cal{R}} \}$$ the pseudo-differential operators given respectively  by
the laws : 
$$v\mapsto (- \Delta)^{1/4}[Qv],~~v\mapsto{\cal{R}}(-\Delta)^{1/4}[Qv],~~ v\mapsto (- \Delta)^{1/4}[Q{\cal{R}}[v]]$$ for $v\in H^{1/2}.$}}
\par
We write 
\begin{equation}\label{MatrixTS}
\left(\begin{array}{c}
 T(P_T,(-\Delta)^{1/4}u)\\ {\cal{R}}[S(P_N,  (-\Delta)^{1/4}u)] 
\end{array}\right)
=C\left( \begin{array}{c} P_T(-\Delta)^{1/4}u\\{\cal{R}}P_N(-\Delta)^{1/4}u\end{array}\right)
\end{equation}
where 
the matrix $C$  is a $2\times 2m$ matrix whose components are made by pseudo-differential operators: for $j\in\{1,\ldots, m\}$
{ \begin{equation}
 \begin{array}{c}
c_{1j}=(-\Delta)^{1/4}\{P_T\}-P_T(-\Delta)^{1/4}+(-\Delta)^{1/4}[P_T]\\[5mm]
c_{1,j+m}=(-\Delta)^{1/4}\{P_T\circ {\cal{R}}\}-P_T(-\Delta)^{1/4}\circ {\cal{R}}+(-\Delta)^{1/4}[P_T]\circ {\cal{R}}\\[5mm]
c_{2,j}={\cal{R}}(-\Delta)^{1/4}\{P_N\}-P_N(-\Delta)^{1/4}\circ {\cal{R}}-(-\Delta)^{1/4}[P_N]\circ {\cal{R}}\\[5mm]
c_{2,j+m}=(-\Delta)^{1/4}\{P_N\}+P_N(-\Delta)^{1/4} -(-\Delta)^{1/4}[P_N]\,. 
\end{array}
\end{equation}}
From \rec{zz7tris} and \rec{zz8tris} it follows that for   $v\in L^2$ the following estimate holds
 
$$
 \|c_{ij}(v)\|_{{{{\cal{H}}}}^{1}(\R)}\le C\ \|P_T\|_{\dot{H}^{1/2}(\R)}\|v\|_{L^{2}(\R)}.
$$}
 {\bf Step 2.} Now we rewrite  the following  matrix
 \begin{equation}\label{matrixF}
\left(\begin{array}{c}
F_1  \\ F_2+F_3   \end{array}\right)
\end{equation}
where
\begin{eqnarray}\label{F}
F_1:=F( \omega_1, P_N \Delta ^{1/4} u)&=&-F(\omega,{\mathcal{R}}[{\mathcal{R}}[P_N \Delta ^{1/4} u]])\\
&=&{\mathcal{R}}[\omega_1]{\mathcal{R}}[P_N(-\Delta)^{1/4}u]+\omega_1{\mathcal{R}}[{\mathcal{R}}[P_N \Delta ^{1/4} u]],\nonumber
 \end{eqnarray}
 
\begin{eqnarray}\label{F2}
F_2 &:=&F(\omega_2,P_N((-\Delta)^{1/4} u)=-F(\omega_2,{\mathcal{R}}[{\mathcal{R}}[P_N \Delta ^{1/4} u]])\\
&=&{\mathcal{R}}[\omega_2]{\mathcal{R}}[P_N(-\Delta)^{1/4}u]+\omega_2{\mathcal{R}}[{\mathcal{R}}[P_N \Delta ^{1/4} u]].\nonumber
\end{eqnarray}
and 
\begin{eqnarray}\label{F3}
&& F_3:=F({\cal{R}}((-\Delta)^{1/4} P_N),{\cal{R}}((-\Delta)^{1/4} u)) \\&=&\underbrace{F({\cal{R}}((-\Delta)^{1/4} P_N),{\cal{R}}P_T((-\Delta)^{1/4} u))}_{(1)}+\underbrace{F({\cal{R}}((-\Delta)^{1/4} P_N),{\cal{R}}P_N((-\Delta)^{1/4} u))}_{(2)}\nonumber
.\end{eqnarray}

We have
 \begin{eqnarray*}
(1)&=& (-\Delta)^{1/4} P_N P_T((-\Delta)^{1/4} u))-{\cal{R}}((-\Delta)^{1/4} P_N) {\cal{R}}P_T((-\Delta)^{1/4} u)\\[5mm]
(2)&=&-((-\Delta)^{1/4} P_N){\cal{R}}({\cal{R}}P_N((-\Delta)^{1/4} u)))-{\cal{R}}((-\Delta)^{1/4} P_N){\cal{R}}P_N((-\Delta)^{1/4} u)).
\end{eqnarray*}
By using  \rec{F}, \rec{F2}, \rec{F3} we can rewrite \rec{matrixF}  in terms of a matrix of pseudo differential operators.

\begin{equation}\label{matrixf}
\left(\begin{array}{c}
F_1  \\ F_2+F_3   \end{array}\right)
=D\left( \begin{array}{c} P_T(-\Delta)^{1/4}u\\{\cal{R}}P_N(-\Delta)^{1/4}u\end{array}\right)
\end{equation}
where $C$ is a matrix of pseudo differential operators given by
{ 
\begin{equation}
 \begin{array}{c}
 d_{1,j}=0\\[5mm]
 d_{1,j+m}={\mathcal{R}}[\omega_1]+\omega_1{\mathcal{R}}\\[5mm]
 d_{2,j}=(-\Delta)^{1/4} P_N+{\cal{R}}((-\Delta)^{1/4} P_N) {\cal{R}}\\[5mm]
 d_{2,j+m}=[{\mathcal{R}}[\omega_2]+\omega_2{\mathcal{R}}]-((-\Delta)^{1/4} P_N){\cal{R}}-{\cal{R}}((-\Delta)^{1/4} P_N).
 \end{array}
 \end{equation}
 Moreover for $v\in L^2$,
 \begin{equation}\label{bij}
\|d_{ij}(v)\|_{  {\cal{H}}^{1}(\R)}\le C \left[\|P_T\|_{\dot H^{1/2}(\R))}+\|P_T\|^2_{\dot H^{1/2}(\R))}\right]\|v\|_{L^{2}}\,.
\end{equation}}

By combining \rec{MatrixTS} and \rec{matrixf} we get
\begin{equation}
{ \tilde\Omega=(C-2D)\left( \begin{array}{c} P_T(-\Delta)^{1/4}u\\{\cal{R}}P_N(-\Delta)^{1/4}u\end{array}\right)}.
\end{equation}
We conclude the proof of Proposition \ref{EulEq}.~~$\Box$

\end{document}